\documentclass{article}
\usepackage{latexsym}
\usepackage{amstext}
\usepackage{amsmath}
\usepackage{amssymb}
\usepackage{amsthm}

\newtheorem{theorem}{Theorem}[section]
\newtheorem{lemma}[theorem]{Lemma}
\newtheorem{theoremb}{Theorem}[section]
\newtheorem{corollary}[theorem]{Corollary}
\newtheorem{corollarya}[theorem]{Corollary}
\newtheorem{corollaryb}[theoremb]{Corollary}
\newtheorem{definition}[theorem]{Definition}
\newtheorem{definitiona}[theorem]{Definition}
\newtheorem{definitionb}[theoremb]{Definition}
\newtheorem{conjecture}[theorem]{Conjecture}
\newtheorem{conjecturea}[theorem]{Conjecture}
\newtheorem{conjectureb}[theoremb]{Conjecture}
\newtheorem{example}[theorem]{Example}
\newtheorem{examplea}[theorem]{Example}
\newtheorem{exampleb}[theoremb]{Example}

\newtheorem{theorema}[theorem]{Theorem}

\theoremstyle{remark}

\newtheorem{remark}[theorem]{Remark}

\renewcommand{\thetheorem}{\arabic{section}.\arabic{theorem}A}

\def \co {\mathcal{O}}
\def \ll {\mathcal{L}}
\def \kk {\overline{k}}
\DeclareMathOperator{\Div}{Div}
\DeclareMathOperator{\ord}{ord}
\DeclareMathOperator{\Sym}{Sym}

\DeclareMathOperator{\Exch}{Exc_{hol}}
\DeclareMathOperator{\Excd}{Exc_{Dio}}
\DeclareMathOperator{\ExcL}{Exc_{Dio,\it{L},S}}
\DeclareMathOperator{\Exck}{Exc_{Dio,\it{k},S}}
\DeclareMathOperator{\dExcd}{Exc_{Dio, deg \it{d}}}
\DeclareMathOperator{\dExcL}{Exc_{Dio,deg \it{d}, \it{L},S}}
\DeclareMathOperator{\dExck}{Exc_{Dio,deg \it{d}, \it{k},S}}

\DeclareMathOperator{\supp}{Supp}

\begin{document}
\bibliographystyle{plain}
\title{Generalizations of Siegel's and Picard's Theorems}
\author{Aaron Levin}
\date{}          
\maketitle
\begin{abstract}
We prove new theorems which are higher-dimensional generalizations of the classical theorems of Siegel on integral points on affine curves and of Picard on holomorphic maps from $\mathbb{C}$ to affine curves.  These include results on integral points over varying number fields of bounded degree and results on Kobayashi hyperbolicity.  We give a number of new conjectures describing, from our point of view, how we expect Siegel's and Picard's theorems to optimally generalize to higher dimensions.  In some special cases we will be able to relate our conjectures to existing conjectures.  In this respect, we are also led to formulate a new conjecture relating the absolute discriminant and height of an algebraic point on a projective variety over a number field.
\end{abstract}
\section{Introduction}
\label{Introduction}
In this article we prove new theorems which are higher-dimensional generalizations of the classical theorems of Siegel on integral points on affine curves and of Picard on holomorphic maps from $\mathbb{C}$ to affine curves.  In the first section  we will give the statements of Siegel's and Picard's theorems, and we will recall why these two theorems from such seemingly different areas of mathematics are related.  We will then proceed to give a number of new conjectures describing, from our point of view, how we expect Siegel's and Picard's theorems to optimally generalize to higher dimensions.  These include conjectures on integral points over varying number fields of bounded degree and conjectures addressing hyperbolic questions.  These conjectures appear to be fundamentally new.  However, in some special cases we will be able to relate our conjectures to Vojta's conjectures.  In this respect, we are also led to formulate a new conjecture relating the absolute discriminant and height of an algebraic point on a projective variety over a number field (Conjecture \ref{conj4}).

We will then summarize our progress on these conjectures.  We have been able to get results in all dimensions, with best-possible results in many cases for surfaces.  Our technique is based on the new proof of Siegel's theorem given by Corvaja and Zannier in \cite{Co}.  They showed how one may use the Schmidt Subspace Theorem to obtain a very simple and elegant proof of Siegel's theorem.  More recently, they have used this technique to obtain other results on integral points  (see \cite{Co5}, \cite{Co3}, and \cite{Co2}) and Ru has translated the approach to Nevanlinna theory \cite{Ru3}.  We will use the Schmidt Subspace Theorem approach to get results on integral points on higher-dimensional varieties, and analogously, we will use a version of Cartan's Second Main Theorem due to Vojta to obtain results on holomorphic curves in higher-dimensional complex varieties, generalizing Picard's theorem.  

As an application of our results, we show how to improve a result of Faltings on integral points on the complements of certain singular plane curves, proving a statement about hyperbolicity as well.  We end with a discussion of our conjectures, relating them to previously known results and conjectures, and giving examples limiting any improvement to their hypotheses and conclusions.

\section{Theorems of Siegel and Picard}
\label{sclassical}
It has been observed by Osgood, Vojta, Lang, and others that there is a striking correspondence between statements in Nevanlinna theory and in Diophantine approximation (see \cite{Ru} and \cite{Vo2}).  This correspondence has been extremely lucrative, influencing results and conjectures in both subjects considerably.  The correspondence can be formulated in both a qualitative and quantitative way.  In this section, we will concentrate on the simplest case of the qualitative correspondence, Siegel's and Picard's theorems.  

Let $V\subset \mathbb{A}^n$ be an affine variety defined over a number field $k$.  We will also view $V$ as a complex analytic space.  Then it has been noticed that $V(\co_{L,S})$ (the set of points with all coordinates in $\co_{L,S}$, the $S$-integers of $L$) seems to be infinite for sufficiently large number fields $L$ and sets of places $S$ if and only if there exists a non-constant holomorphic map $f:\mathbb{C}\to V$.  When $V=C$ is a curve (i.e. one-dimensional variety), this correspondence has been proven to hold exactly, and it is known precisely for which curves $C$ the two statements hold.  On the number theory side, Siegel's theorem is the fundamental theorem on integral points on curves.  On the analytic side the analogue is a theorem of Picard.  We now give the following formulations of these two theorems.

\begin{theorema}[Siegel]
\label{Siegel2}
Let $k$ be a number field.  Let $S$ be a finite set of places of $k$ containing the archimedean places.  Let $C$ be an affine curve defined over $k$ embedded in affine space $\mathbb{A}^m$.  Let $\tilde{C}$ be a projective closure of $C$.  If $\# \tilde{C}\backslash C >2$ (over $\kk$) then $C$ has finitely many points in $\mathbb{A}^m(\mathcal{O}_{k,S})$.
\end{theorema}
\begin{theoremb}[Picard]
\label{Picard}
Let $\tilde{C}$ be a compact Riemann surface.  Let $C \subset \tilde{C}$.  If $\# \tilde{C}\backslash C > 2$, then all holomorphic maps $f:\mathbb{C} \to C$ are constant.
\end{theoremb}

In other words, Siegel's and Picard's theorems state that if $D$ consists of many distinct points on a curve $X$, then any set of integral points on $X\backslash D$ is finite and any holomorphic map $f:\mathbb{C}\to X\backslash D$ is constant.  We will thus view as generalizing Siegel's or Picard's theorem any theorem that asserts that if $D$ has ``enough components" then there is some limitation on the integral points on $X\backslash D$ or on the holomorphic maps $f:\mathbb{C}\to X\backslash D$.  In Picard's theorem it may also be shown that the curves $C$ in question satisfy the stronger condition of being Kobayashi hyperbolic.  We will frequently be able to generalize this fact to higher dimensions as well.

Siegel's theorem is usually stated with the extra information that the $\# \tilde{C}\backslash C >2$ hypothesis is unnecessary for nonrational affine curves $C$.  However, it may be shown that this stronger version of Siegel's theorem may be derived from Siegel's theorem as we have stated it by using \'etale coverings of the curve $C$ (see \cite{Co}).  A similar statement holds for Picard's theorem.  It is Siegel's and Picard's theorems in the form we have given above that we will generalize.

We note that when the geometric genus of $C$ is greater than one, Siegel's theorem follows from the much stronger theorem of Faltings that $C$ has only finitely many $k$-rational points.  Similarly, it is a theorem of Picard that there are no nonconstant holomorphic maps $f:\mathbb{C}\to \tilde{C}$ when $\tilde{C}$ is a projective curve of geometric genus greater than one.  
\section{Some Preliminary Definitions}
In order to state our conjectures and results we will need a few definitions.  In Vojta's Nevanlinna-Diophantine dictionary, the Diophantine object corresponding to a holomorphic map $f:\mathbb{C}\to X\backslash D$ is a set of $(D,S)$-integral points on $X$.  We'll now sketch the definition of a set of $(D,S)$-integral points on $X$ in terms of Weil functions.

Let $D$ be a Cartier divisor on a projective variety $X$, both defined over a number field $k$.  Let $M_k$ denote the set of places of $k$ (see Section \ref{sDio}).  Let $v\in M_k$.  Extend $|\cdot|_v$ to an absolute value on $\kk_v$.  We define a local Weil function for $D$ relative to $v$ to be a function $\lambda_{D,v}:X(\kk_v)\backslash D \to \mathbb{R}$ such that if $D$ is represented locally by $(f)$ on an open set $U$ then
\renewcommand{\thetheorem}{\arabic{section}.\arabic{theorem}}
\begin{equation*}
\lambda_{D,v}(P)=-\log{|f(P)|_v}+\alpha_v(P)
\end{equation*}
where $\alpha_v$ is a continuous function on $U(\kk_v)$ (in the $v$-topology).  

By choosing embeddings $k\to \kk_v$ and $\kk\to \kk_v$, we may also think of $\lambda_{D,v}$ as a function on $X(k)\backslash D$ or $X(\kk)\backslash D$.  A global Weil functions consists of a collection of local Weil functions, $\lambda_{D,v}$, for $v\in M_k$, where the $\alpha_v$ above satisfy certain reasonable boundedness conditions as $v$ varies.  We refer the reader to \cite{La} and \cite{Vo2} for a further discussion of this.
\begin{definition}
Let $D$ be an effective Cartier divisor on a projective variety $X$, both defined over a number field $k$.  Let $S$ be a finite set of places in $M_k$ containing the archimedean places.  Let $R \subset X(\kk)\backslash D$.  Then $R$ is defined to be a $(D,S)$-integral set of points if there exists a global Weil function $\lambda_{D,v}$ and constants $c_v$, with $c_v=0$ for all but finitely many $v$, such that for all $v\in M_k\backslash S$ and all embeddings $\kk \to \kk_v$\\
\begin{equation*}
\lambda_{D,v}(P) \leq c_v
\end{equation*}
for all $P$ in $R$.
\end{definition}
We will frequently just say $D$-integral, omitting the reference to $S$, when $S$ has been fixed or when the statement is true for all possible $S$.  Except where explicitly stated, we will also require from now on that a set of $D$-integral points be $k$-rational, i.e. $R\subset X(k)$.

For us, the key property of a set of $(D,S)$-integral points is given by the following theorem.
\begin{theorem}
\label{reg}
Let $R\subset X(\kk)\backslash D$ be a set of $(D,S)$-integral points on $X$.  Then for any regular function $f$ on $X\backslash D$ (defined over $\kk$) there exists a constant $a\in k$ such that $af(P)$ is $S$-integral for all $P$ in $R$, that is $af(P)$ lies in the integral closure of $\co_{k,S}$ in $\kk$ for all $P\in R$.
\end{theorem}
In fact, in what follows, most of our results hold, and our conjectures should hold, for any set $R$ satisfying the conclusion of Theorem \ref{reg}.  We will prefer to work with sets of $D$-integral points because they are better geometrically behaved (e.g. under pullbacks) and because they are the right objects to use so that the Diophantine exceptional set we are about to define matches (conjecturally) the holomorphic exceptional set we will define.  We note that sets of $D$-integral points are also essentially the same as the sets of scheme-theoretic integral points one would get from working with models of $X\backslash D$ over $\co_{k,S}$ (see \cite[Prop. 1.4.1]{Vo2}).

It will be necessary to define various exceptional sets of a variety.
\renewcommand{\thetheorem}{\arabic{section}.\arabic{theorem}A}
\setcounter{theoremb}{\value{theorem}}
\begin{definitiona}
Let $X$ be a projective variety and $D$ an effective Cartier divisor on $X$, both defined over a number field $k$.  Let $L$ be a number field, $L\supset k$, and $S$ a finite set of places of $L$ containing the archimedean places.  We define the Diophantine exceptional set of $X\backslash D$ with respect to $L$ and $S$ to be
\begin{equation*}
\ExcL(X\backslash D)=\bigcup_R \dim_{>0}(\overline{R})
\end{equation*}
 where the union runs over all sets $R$ of $L$-rational $(D,S)$-integral points on $X$ and $\dim_{>0}(\overline{R})$ denotes the union of the positive dimensional irreducible components of the Zariski-closure of $R$.  We define the absolute Diophantine exceptional set of $X\backslash D$ to be 
\begin{equation*}
\Excd(X\backslash D)=\bigcup_{L \supset k,S} \ExcL(X\backslash D),
\end{equation*}
with $L$ ranging over all number fields and $S$ ranging over all sets of places of $L$ as above.
\end{definitiona}
These definitions depend only on $X\backslash D$ and not on the choices of $X$ and $D$.
\begin{definitionb}
Let $X$ be a complex variety.  We define the holomorphic exceptional set of $X$, $\Exch(X)$, to be the union of all images of non-constant holomorphic maps $f:\mathbb{C}\to X$.
\end{definitionb}
Conjecturally, it is expected that $\Excd(X\backslash D)=\Exch(X\backslash D)$ (it may also be necessary to take the Zariski-closures of both sides first).  
\begin{definitiona}
Let $X$ be a projective variety defined over a number field $k$.  Let $D$ be an effective Cartier divisor on $X$.  Then we define $X\backslash D$ to be Mordellic if $\Excd(X\backslash D)$ is empty.  We define $X\backslash D$ to be quasi-Mordellic if $\Excd(X\backslash D)$ is not Zariski-dense in $X$. 
\end{definitiona}
\begin{definitionb}
Let $X$ be a complex variety.  We define $X$ to be Brody hyperbolic if $\Exch(X)$ is empty.  We define $X$ to be quasi-Brody hyperbolic if $\Exch(X)$ is not Zariski-dense in $X$.
\end{definitionb}
Note that $X$ being quasi-Brody hyperbolic is a stronger condition than the non-existence of holomorphic maps $f:\mathbb{C}\to X$ with Zariski-dense image.  Similarly, $X\backslash D$ being quasi-Mordellic is stronger than the non-existence of dense sets of $D$-integral points on $X$.

We will also need a convenient measure of the size of a divisor.  We will use $\co_X(D)$, or simply $\co(D)$ when there is no ambiguity, to denote the invertible sheaf associated to a Cartier divisor $D$ on $X$, and $h^i(D)$ to denote the dimension of the vector space $H^i(X,\co(D))$.  When $h^0(D)>0$, we will also frequently use the notation $\Phi_D$ to denote the rational map (unique up to projective automorphisms) from $X$ to $\mathbb{P}^{h^0(D)-1}$ corresponding to a basis of $H^0(X,\co(D))$.
\renewcommand{\thetheorem}{\arabic{section}.\arabic{theorem}}
\begin{definition}
\label{defk}
Let $D$ be a divisor on a nonsingular projective variety $X$.  We define the dimension of $D$ to be the integer $\kappa(D)$ such that there exists positive constants $c_1$ and $c_2$ such that
\begin{equation*}
c_1 n^{\kappa(D)} \leq h^0(nD)\leq c_2 n^{\kappa(D)}
\end{equation*}
for all sufficiently divisible $n>0$.  If $h^0(nD)=0$ for all $n>0$ then we let $\kappa(D)=-\infty$.
\end{definition}
Alternatively, if $\kappa(D)\geq 0$, one can show that 
\begin{equation*}
\kappa(D)=\max \{\dim \Phi_{nD}(X)|n>0,h^0(nD)>0\}.
\end{equation*}
If $D$ is a Cartier divisor on a singular complex projective variety, we define $\kappa(D)=\kappa(\pi^*D)$ where $\pi:X'\to X$ is a desingularization of $X$.  It is easy to show that this is independent of the chosen desingularization.  For more properties of $\kappa(D)$ we refer the reader to \cite[Ch. 10]{Ii}.
\begin{definition}
\label{defbig}
We define a Cartier divisor $D$ on $X$ to be quasi-ample (or big) if $\kappa(D)=\dim X$.
\end{definition}
If $D$ is quasi-ample then there exists an $n>0$ such that $\Phi_{nD}$ is birational, justifying the name.

\section{General Setup and Notation}
\label{gsetup}
Throughout this paper we will use the following general setup and notation.\\\\
\textbf {General setup}:  Let $X$ be a complex projective variety.  Let $D=\sum_{i=1}^r D_i$ be a divisor on $X$ with the $D_i$'s effective Cartier divisors for all $i$.  Suppose that at most $m$ $D_i$'s meet at a point, so that the intersection of any $m+1$ distinct $D_i$'s is empty.\\\\
In the Diophantine setting, we will also assume that $X$ and $D$ are defined over a number field $k$ and we let $S$ be a finite set of places of $k$ containing the archimedean places.

From now on, we will freely use the notation $X$, $D$, $D_i$, $r$, $m$, $k$, and $S$ as above without further explanation.

\section{Siegel and Picard-type Conjectures}
In this section we give conjectures generalizing Siegel's theorem and Picard's theorem in various directions.
\subsection{Main Conjectures}
Some special cases of the conjectures given in this section are related to Vojta's Main Conjecture.  Later, we will also give conjectures related to Vojta's General Conjecture, hence our terminology in this section and the next.  We remind the reader that throughout we are using the general setup of the last section.
\renewcommand{\thetheorem}{\arabic{section}.\arabic{theorem}A}
\setcounter{theoremb}{\value{theorem}}
\begin{conjecturea}[Main Siegel-type Conjecture]
\label{conjmaina}
Suppose that $\kappa(D_i)\geq \kappa_0>0$ for all $i$.  If $r>m+\frac{m}{\kappa_0}$ then there does not exist a Zariski-dense set of $k$-rational $(D,S)$-integral points on $X$.
\end{conjecturea}
\begin{conjectureb}[Main Picard-type Conjecture]
\label{conjmainb}
Suppose that $\kappa(D_i)\geq \kappa_0>0$ for all $i$.  If $r>m+\frac{m}{\kappa_0}$ then there does not exist a holomorphic map $f:\mathbb{C} \to X \backslash D$ with Zariski-dense image.
\end{conjectureb}
As mentioned earlier, we will usually just say $D$-integral, omitting $k$ and $S$ from the notation.
Siegel's theorem (resp. Picard's theorem) is the case $m=\kappa_0=\dim X=1$ of Conjecture \ref{conjmaina} (resp. Conjecture \ref{conjmainb}).  We note that the dimension of $X$ does not appear in the conjectures, but $\kappa(D_i)$ is bounded by $\dim X$.  We will now discuss some consequences and special cases of these conjectures which seem important enough in their own right to be listed separately as new conjectures, and which will sometimes contain extra conjectures (e.g. on the exceptional sets) which do not follow from the main conjectures above.  At the two extremes of $\kappa_0$ we have
\begin{conjecturea}
\label{conj1a}
If $\kappa(D_i)>0$ for all $i$ and $r> 2m$ then there does not exist a Zariski-dense set of $D$-integral points on $X$.
\end{conjecturea}
\begin{conjectureb}
\label{conj1b}
If $\kappa(D_i)>0$ for all $i$ and $r>2m$ then there does not exist a holomorphic map $f:\mathbb{C} \to X \backslash D$ with Zariski-dense image.
\end{conjectureb}
\begin{conjecturea}
\label{conj1ab}
If $D_i$ is quasi-ample for all $i$ and $r>m+\frac{m}{\dim X}$ then $X\backslash D$ is quasi-Mordellic.
\end{conjecturea}
\begin{conjectureb}
\label{conj1bb}
If $D_i$ is quasi-ample for all $i$ and $r>m+\frac{m}{\dim X}$ then $X\backslash D$ is quasi-Brody hyperbolic.
\end{conjectureb}
We note that when the $D_i$'s are in some sort of general position, so that $m=\dim X$, the inequalities in the last two conjectures above take the nicer form $r>\dim X +1$.  The statements on quasi-Mordellicity and quasi-Brody hyperbolicity do not follow (directly at least) from the Main Conjectures.

Of particular interest is the case where $D_i$ is ample for all $i$.  In this case we conjecture very precise bounds on the dimensions of the exceptional sets (see Remark \ref{rbig} for a possible generalization to quasi-ample divisors).
\begin{conjecturea}[Main Siegel-type Conjecture for Ample Divisors]
\label{conj2a}
Suppose that $D_i$ is ample for all $i$.\\\\
(a).  If $r>m+\frac{m}{\dim X}$ then $\dim \Excd(X) \leq \frac{m}{r-m}$.\\
(b).  In particular, if $r>2m$ then $X\backslash D$ is Mordellic.
\end{conjecturea}
\begin{conjectureb}[Main Picard-type Conjecture for Ample Divisors]
\label{conj2b}
Suppose that $D_i$ is ample for all $i$.\\\\
(a).  If $r>m+\frac{m}{\dim X}$ then $\dim \Exch(X) \leq \frac{m}{r-m}$.\\
(b).  If $r>2m$ then $X\backslash D$ is complete hyperbolic and hyperbolically imbedded in $X$.  In particular, $X\backslash D$ is Brody hyperbolic.
\end{conjectureb}
It is not hard to show that the Main Conjectures for ample divisors follow from Conjectures \ref{conj1ab} and \ref{conj1bb}.

\subsection{General Conjectures}
We will also consider the situation where the field that the integral points are defined over is allowed to vary over all fields of degree less than or equal to $d$ over some fixed field $k$.  So in this section we do not require that the integral points be $k$-rational.
\renewcommand{\thetheorem}{\arabic{section}.\arabic{theorem}}
\begin{definition}
Let $R\subset X(\kk)$.  We define the degree of $R$ over $k$ to be $\deg_k R=\sup_{P\in R} [k(P):k]$.
\end{definition}

Generalizing the Main Siegel-type Conjecture of last section, we conjecture
\begin{conjecture}[General Siegel-type Conjecture]
\label{congen}
Suppose that $\kappa(D_i)\geq \kappa_0>0$ for all $i$.  Let $d$ be a positive integer.  If $r>m+\frac{m(2d-1)}{\kappa_0}$ then there does not exist a Zariski-dense set of $D$-integral points on $X$ of degree $d$ over $k$.
\end{conjecture}
\noindent We will see later that this conjecture and others in this section are related to Vojta's General Conjecture.

We will also want to define a degree $d$ Diophantine exceptional set for a variety $V$.  With the notation from our earlier definition for $\Excd$ we define
\begin{definition}
Let $X$ be a projective variety and $D$ an effective Cartier divisor on $X$, both defined over a number field $k$.  Let $L$ be a number field, $L\supset k$, and $S$ a finite set of places of $L$ containing the archimedean places.  We define the degree $d$ Diophantine exceptional set of $X\backslash D$ with respect to $L$ and $S$ to be
\begin{equation*}
\dExcL(X\backslash D)=\bigcup_R \dim_{>0}(\overline{R})
\end{equation*}
 where the union runs over all sets $R$ of $(D,S)$-integral points on $X$ of degree $d$ over $L$.  We define the degree $d$ absolute Diophantine exceptional set of $X\backslash D$ to be 
\begin{equation*}
\dExcd(X\backslash D)=\bigcup_{L \supset k,S} \dExcL(X\backslash D),
\end{equation*}
with $L$ ranging over all number fields and $S$ ranging over all sets of places of $L$ as above.
\end{definition}
Similarly we define $X\backslash D$ to be degree $d$ Mordellic (resp. degree $d$ quasi-Mordellic) if $\dExcd(X\backslash D)$ is empty (resp. not Zariski-dense in X).  At the two extremes of $\kappa_0$ we have
\begin{conjecture}
Let $d$ be a positive integer.  If $\kappa(D_i)>0$ for all $i$ and $r>2dm$ then there does not exist a Zariski-dense set of $D$-integral points on $X$ of degree $d$ over $k$.
\end{conjecture}
\begin{conjecture}
Let $d$ be a positive integer.  If $D_i$ is quasi-ample for all $i$ and $r>m+\frac{m(2d-1)}{\dim X}$ then $X\backslash D$ is degree $d$ quasi-Mordellic.
\end{conjecture}
We can also give a conjecture for ample divisors giving bounds on the degree $d$ Diophantine exceptional set.
\begin{conjecture}[General Siegel-type Conjecture for Ample Divisors]
Suppose that $D_i$ is ample for all $i$.\\\\
(a).  If $r>m+\frac{m(2d-1)}{\dim X}$ then $\dim \dExcd(X\backslash D)\leq \frac{m(2d-1)}{r-m}$.\\
(b).  In particular, if $r>2dm$ then $X\backslash D$ is degree $d$ Mordellic.
\end{conjecture}

\subsection{Conjectures over $\mathbb{Z}$ and Complex Quadratic Rings of Integers}
When $\#S=1$, or equivalently, when $\co_{k,S}$ is $\mathbb{Z}$ or the ring of integers of a complex quadratic field, and $D_i$ is defined over $k$ for all $i$, we conjecture improvements to our previous conjectures.  We will refer to these conjectures as ``over  $\mathbb{Z}$", though they apply equally well to rings of integers of complex quadratic fields.
\begin{conjecture}[Main Siegel-type Conjecture over $\mathbb{Z}$]
Let $k=\mathbb{Q}$ or a complex quadratic field and let $S=\{v_\infty\}$ consist of the unique archimedean place of $k$.  Suppose that $D_i$ is defined over $k$ for all $i$ and that $\kappa(D_i)>0$ for all $i$.  If $r>m$ then there does not exist a Zariski-dense set of $(D,S)$-integral points on $X$.
\end{conjecture}
We emphasize that in contrast to our previous conjectures, each $D_i$ must be defined over $k$.  We also conjecture that in the above if each $D_i$ is quasi-ample, then $\Exck(X)$ is not Zariski-dense in $X$.  For ample divisors, as usual, we conjecture something more.
\begin{conjecture}[Main Siegel-type Conjecture over $\mathbb{Z}$ for Ample Divisors]
\label{conjS}
Let $k=\mathbb{Q}$ or a complex quadratic field and let $S=\{v_\infty\}$ consist of the unique archimedean place of $k$.  Suppose that $D_i$ is ample and defined over $k$ for all $i$.\\\\
(a).  All sets $R$ of $(D,S)$-integral points on $X$ have $\dim R \leq 1+\dim (\bigcap_i D_i)$.\\
(b).  In particular, if $D=D_1+D_2$ is a sum of two ample effective Cartier divisors on $X$, both defined over $k$, with no irreducible components in common,  then there does not exist a Zariski-dense set of $(D,S)$-integral points on $X$. 
\end{conjecture}
\begin{conjecture}[General Siegel-type Conjecture over $\mathbb{Z}$]
\label{GZ}
Let $k=\mathbb{Q}$ or a complex quadratic field and let $S=\{v_\infty\}$ consist of the unique archimedean place of $k$.  Suppose that $D_i$ is defined over $k$ for all $i$ and that $\kappa(D_i)\geq \kappa_0>0$ for all $i$.  Let $d$ be a positive integer.  If $r>m+\frac{m(d-1)}{\kappa_0}$ then there does not exist a Zariski-dense set of $(D,S)$-integral points on $X$ of degree $d$ over $k$.
\end{conjecture}
If $D_i$ is quasi-ample for all $i$ in the above conjecture, then we also conjecture that $\dExck(X\backslash D)$ is not Zariski-dense in $X$.  For ample divisors we have
\begin{conjecture}[General Siegel-type Conjecture over $\mathbb{Z}$ for Ample Divisors]
Let $k=\mathbb{Q}$ or a complex quadratic field and let $S=\{v_\infty\}$ consist of the unique archimedean place of $k$.  Suppose that $D_i$ is ample and defined over $k$ for all $i$.  Let $d$ be a positive integer.  If  $r>m+\frac{m(d-1)}{\dim X}$ then $\dim \dExck(X\backslash D)\leq \frac{m(d-1)}{r-m}$.
\end{conjecture}
\renewcommand{\thetheorem}{\arabic{section}.\arabic{theorem}A}
\setcounter{theoremb}{\value{theorem}}
We will discuss the conjectures in greater detail in Section \ref{Remarks}.
\section{Overview of Results}
Sections \ref{smain}-\ref{SVgeneral} will be concerned with proving special cases of the above conjectures.  In this section we highlight some of our results.  
Along the lines of the Main Conjectures we have
\begin{theorema}
Suppose $r>2m\dim X$.\\\\
(a).  If $D_i$ is quasi-ample for all $i$ then $X\backslash D$ is quasi-Mordellic.\\
(b).  If $D_i$ is ample for all $i$ then $X\backslash D$ is Mordellic.
\end{theorema}
\begin{theoremb}
Suppose $r>2m\dim X$.\\\\
(a).  If $D_i$ is quasi-ample for all $i$ then $X\backslash D$ is quasi-Brody hyperbolic.\\
(b).  If $D_i$ is ample for all $i$ then $X\backslash D$ is complete hyperbolic and hyperbolically imbedded in $X$.  In particular, $X\backslash D$ is Brody hyperbolic.
\end{theoremb}
If in addition the $D_i$'s have no irreducible components in common, then in the part (a)'s above we only need  $r>2[\frac{m+1}{2}]\dim X$ where $[x]$ denotes the greatest integer in $x$.

When $X$ is a surface, $m\leq 2$, and the $D_i$'s have no irreducible components in common, we are able to prove the Main Conjectures, Conjectures \ref{conjmaina},B through \ref{conj2a},B.
\begin{theorema}
Suppose $X$ is a surface and the $D_i$'s have no irreducible components in common.\\\\
(a).  If $m=1$, $\kappa(D_i)>0$ for all $i$, and $r>2$ then there does not exist a Zariski-dense set of $D$-integral points on $X$.\\
(b).  If $m=2$, $\kappa(D_i)>0$ for all $i$, and $r>4$ then there does not exist a Zariski-dense set of $D$-integral points on $X$.\\
(c).  If $m=2$, $D_i$ is quasi-ample for all $i$, and $r>3$ then $X\backslash D$ is quasi-Mordellic.\\
(d).  If $m=2$, $D_i$ is ample for all $i$, and $r>4$ then $X\backslash D$ is Mordellic.  
\end{theorema}
\begin{theoremb}
Suppose $X$ is a surface and the $D_i$'s have no irreducible components in common.\\\\
(a).  If $m=1$, $\kappa(D_i)>0$ for all $i$, and $r>2$ then there does not exist a holomorphic map $f:\mathbb{C}\to X\backslash D$ with Zariski-dense image.\\
(b).  If $m=2$, $\kappa(D_i)>0$ for all $i$, and $r>4$ then there does not exist a holomorphic map $f:\mathbb{C}\to X\backslash D$ with Zariski-dense image.\\
(c).  If $m=2$, $D_i$ is quasi-ample for all $i$, and $r>3$ then $X\backslash D$ is quasi-Brody hyperbolic.\\
(d).  If $m=2$, $D_i$ is ample for all $i$, and $r>4$ then $X\backslash D$ is complete hyperbolic and hyperbolically imbedded in $X$.  In particular, $X\backslash D$ is Brody hyperbolic.
\end{theoremb}
We will see later that if $m=1, r>1,$ and $\kappa(D_i)>0$ for all $i$, then we must necessarily have $\kappa(D_i)=1$ for all $i$.

As to the General Conjectures, when the integral points are allowed to vary over fields of a bounded degree, we prove
\renewcommand{\thetheorem}{\arabic{section}.\arabic{theorem}}
\begin{theorem}
Let $d$ be a postive integer.  If $D_i$ is ample for all $i$ and $r>2d^2m\dim X$ then $X\backslash D$ is degree $d$ Mordellic (all sets of $D$-integral points on $X$ of degree $d$ over $k$ are finite).
\end{theorem}
\begin{theorem}
Let $k=\mathbb{Q}$ or a complex quadratic field.  Let $S=\{v_\infty\}$ consist of the unique archimedean place of $k$.  Let $d$ be a positive integer.  If $D_i$ is ample and defined over $k$ for all $i$ and $r>dm$ then all sets of $(D,S)$-integral points on $X$ of degree $d$ over $k$ are finite.
\end{theorem}

As an application of our results, we will discuss an improvement to a result of Faltings.
Faltings \cite{Fa} has recently shown how theorems on integral points on the complements of divisors with many components may occasionally be used to prove theorems on the complements of irreducible divisors.  He shows how to do this with certain very singular curves on $\mathbb{P}^2$ by reducing the problem to a covering surface and applying the method of \cite{Fa2}.  In \cite{Co4}, Zannier uses the subspace theorem approach instead of \cite{Fa2} to prove a result similar to Faltings.  In Section \ref{Faltings} we will prove a theorem which generalizes both results.  As an added bonus, we also prove the theorem in the case of holomorphic curves.

\renewcommand{\thetheorem}{\arabic{section}.\arabic{theorem}A}
\setcounter{theoremb}{\value{theorem}}
\section{Preliminaries}
\subsection{Diophantine Approximation}
\label{sDio}
\renewcommand{\thetheorem}{\arabic{section}.\arabic{theorem}}
Let $k$ be a number field.  Let $\mathcal{O}_k$ be the ring of integers of $k$.  As usual, we have a set $M_k$ of absolute values (or places) of $k$ consisting of one place for each prime ideal $\mathfrak{p}$ of $\mathcal{O}_k$, one place for each real embedding $\sigma:k \to \mathbb{R}$, and one place for each pair of conjugate embeddings $\sigma,\overline{\sigma}:k \to \mathbb{C}$.  Let $k_v$ denote the completion of $k$ with respect to $v$.  We normalize our absolute values so that $|p|_v=p^{-[k_v:\mathbb{Q}_p]/[k:\mathbb{Q}]}$ if $v$ corresponds to $\mathfrak{p}$ and $\mathfrak{p}|p$, and $|x|_v=|\sigma(x)|^{[k_v:\mathbb{R}]/[k:\mathbb{Q}]}$ if $v$ corresponds to an embedding $\sigma$ (in which case we say that $v$ is archimedean).  If $v$ is a place of $k$ and $w$ is a place of a field extension $L$ of $k$, then we say that $w$ lies above $v$, or $w|v$, if $w$ and $v$ define the same topology on $k$.

With the above definitions we have the product formula
\begin{equation*}
\prod_{v \in M_k}|x|_v=1 \quad \text{for all } x\in k^*.
\end{equation*}
For a point $P=(x_0,\ldots,x_n)\in \mathbb{P}^n(k)$ we define the height to be
\begin{equation*}
H(P)=\prod_{v\in M_k} \max(|x_0|_v,\ldots,|x_n|_v).
\end{equation*}
It follows from the product formula that $H(P)$ is independent of the choice of homogeneous coordinates for $P$.  It is also easy to see that the height is independent of $k$.  We define the logarithmic height to be
\begin{equation*}
h(P)=\log H(P).
\end{equation*}

At the core of our Diophantine results is the following version of Schmidt's Subspace Theorem due to Vojta \cite{Vo6}.
\renewcommand{\thetheorem}{\arabic{section}.\arabic{theorem}A}
\setcounter{theoremb}{\value{theorem}}
\begin{theorema}
Let $k$ be a number field.  Let $S$ be a finite set of places in $M_k$ containing the archimedean places.  Let $H_1,\ldots,H_m$ be hyperplanes in $\mathbb{P}^n$ defined over $\kk$ with corresponding Weil functions $\lambda_{H_1},\ldots,\lambda_{H_m}$.  Then there exists a finite union of hyperplanes $Z$, depending only on $H_1,\ldots,H_m$ (and not $k$ or $S$),  such that for any $\epsilon>0$,
\begin{equation}
\sum_{v\in S}\max_I \sum_{i \in I} \lambda_{H_i,v}(P) \leq (n+1+\epsilon)h(P)
\end{equation}
holds for all but finitely many $P$ in $\mathbb{P}^n(k)\backslash Z$, where the max is taken over subsets $I \subset \{1,\ldots,m\}$ such that the linear forms defining $H_i,i \in I$ are linearly independent.
\end{theorema}
Explicitly, if $H$ is a hyperplane on $\mathbb{P}^n$ defined by the linear form $L(x_0,\ldots,x_n)$ then a Weil function for $H$ is given by
\begin{equation}
\label{Weila}
\lambda_{H,v}(P)=\log \max_i \frac{|x_i|_v}{|L(P)|_v}.
\end{equation}
where $P=(x_0,\cdots,x_n)$.

We will also need the close relative of Schmidt's theorem, the $S$-unit lemma.
\begin{theorema}
Let $k$ be a number field and let $n>1$ be an integer.  Let $\Gamma$ be a finitely generated subgroup of $k^*$.  Then all but finitely many solutions of the equation
\begin{equation}
u_0+u_1+\cdots+u_n=1, u_i\in \Gamma
\end{equation}
lie in one of the diagonal hyperplanes $H_I$ defined by the equation $\sum_{x\in I}x_i=0$, where $I$ is a proper subset of $\{0,\ldots,n\}$ with at least two elements.
\end{theorema}

\renewcommand{\thetheorem}{\arabic{section}.\arabic{theorem}}
For the convenience of the reader, we have collected various properties of $D$-integral points that we will use (sometimes implicitly) throughout the paper (see \cite{Vo2}).
\begin{lemma}
\label{Integral}
Let $k$ be a number field and $S$ a finite set of places in $M_k$ containing the archimedean places.  Let $D$ be an effective Cartier divisor on a projective variety $X$, both defined over $k$.\\\\
(a).  Let $L$ be a finite extension of $k$ and let $T$ be the set of places of $L$ lying over places in $S$.  If $R$ is a set of $(D,S)$-integral points then it is a set of $(D,T)$-integral points.\\
(b).  Let $E$ be an effective Cartier divisor on $X$.  If $R$ is a set of $(D+E)$-integral points then $R$ is a set of $D$-integral points.\\
(c).  The $D$-integrality of a set is independent of the multiplicities of the components of $D$.\\
(d).  Let $Y$ be a projective variety defined over $k$.  Let $\pi:Y\to X$ be a morphism defined over $k$ with $\pi(Y)\not\subset D$.  If $R$ is a set of $(D,S)$-integral points on $X$ then $\pi^{-1}(R)$ is a set of $(\pi^*D,S)$-integral points on $Y$.
\end{lemma}
Note also in (d), that if in addition $\pi:Y\backslash \pi^*D \to X \backslash D$ is a finite \'etale map, then by the Chevally-Weil theorem there exists a number field $L$ such that $\pi^{-1}(R)\subset Y(L)$ \cite[Th. 1.4.11]{Vo2}.

\subsection{Nevanlinna Theory and Kobayashi Hyperbolicity}
We will be interested in Nevanlinna theory as it applies to holomorphic maps $f:\mathbb{C} \to \mathbb{P}^n$ and hyperplanes on $\mathbb{P}^n$.
Let $f:\mathbb{C} \to \mathbb{P}^n$ be a holomorphic map.  Then we may choose a representation of $f$, $\mathbf{f}=(f_0,\ldots,f_n)$ where $f_0,\ldots,f_n$ are entire functions without common zeros.  Let us define $\|\mathbf{f}\|=(|f_0|^2+\cdots +|f_n|^2)^{\frac{1}{2}}$.  Then we define a characteristic function $T_f(r)$ of $f$ to be
\begin{equation*}
T_f(r)=\int_{0}^{2\pi} \log \|\mathbf{f}(re^{i\theta})\|\frac{d\theta}{2\pi}.
\end{equation*}
Note that by Jensen's formula this function is well-defined up to a constant.  Let $H$ be a hyperplane in $\mathbb{P}^n$ defined by a linear form $L$.  Then we define a Weil function $\lambda_H(f(z))$ of $f$ with respect to $H$ by
\begin{equation}
\label{Weilb}
\lambda_H(f(z))=-\log \frac{|L(\mathbf{f}(z))|}{\|\mathbf{f}(z)\|}.
\end{equation}
We note that this is independent of the choice of $\mathbf{f}$ and depends on the choice of $L$ only up to a constant.  The analogue of Schmidt's Subspace Theorem that we will need is the following version of Cartan's Second Main Theorem, due to Vojta \cite{Vo3}.
\begin{theoremb}
Let $H_1,\ldots H_m$ be hyperplanes in $\mathbb{P}^n$ with corresponding Weil functions $\lambda_{H_1},\ldots,\lambda_{H_m}$.  Then there exists a finite union of hyperplanes $Z$ such that for any $\epsilon >0$ and any holomorphic map $f:\mathbb{C}\to \mathbb{P}^n\backslash Z$
\begin{equation}
\int_{0}^{2\pi} \max_I \sum_{i \in I} \lambda_{H_i}(f(re^{i\theta}))\frac{d\theta}{2\pi} \leq (n+1+\epsilon)T_f(r)
\end{equation}
holds for all $r$ outside a set of finite Lebesgue measure, where the max is taken over subsets $I \subset \{1,\ldots,m\}$ such that the linear forms defining $H_i,i \in I$, are linearly independent.
\end{theoremb}
The analogue of the $S$-unit lemma is the Borel lemma.
\begin{theoremb}
Let $f_1,\ldots,f_n$ be entire functions.  Suppose that
\begin{equation}
e^{f_1}+\cdots+e^{f_n}=1.
\end{equation}
Then $f_i$ is constant for some $i$.
\end{theoremb}

Closely connected to questions about holomorphic curves is the Kobayashi pseudo-distance and Kobayashi hyperbolicity.  We refer the reader to \cite{La2} for the definitions of the Kobayashi pseudo-distance, Kobayashi hyperbolic, complete hyperbolic, and hyperbolically imbedded.  It is trivial that Kobayashi hyperbolic implies Brody hyperbolic.  We will want a criterion for proving the converse in special cases.  On projective varieties, this is given by Brody's theorem.  More generally, we will use the following theorem of Green (see \cite{Gr2} and \cite{La2}).
\renewcommand{\thetheorem}{\arabic{section}.\arabic{theorem}}
\begin{theorem}[Green]
\label{hyperbolic}
Let $X$ be a complex projective variety.  Let $Y=\bigcup_{i \in I} D_i$ be a finite union of Cartier divisors $D_i$ on $X$.  Suppose that for every subset $\emptyset \subset J \subset I$,
\begin{equation*}
\bigcap_{j\in J}D_j\backslash \bigcup_{i\in I\backslash J}D_i
\end{equation*}
is Brody hyperbolic, where $\bigcap_{j\in \emptyset}D_j=X$.  Then $X\backslash Y$ is complete hyperbolic and hyperbolically imbedded in $X$.
\end{theorem}
\subsection{Nef and Quasi-ample Divisors}
\renewcommand{\thetheorem}{\arabic{section}.\arabic{theorem}}
We now recall some basic definitions and facts regarding nef and  quasi-ample divisors.  We will use the theory of intersection numbers on projective varieties as presented in, for instance, \cite{Kl}.  We will use the notation $D^n$ to denote the intersection number of the $n$-fold intersection of $D$ with itself.  In what follows $X$ will be a projective variety over an algebraically closed field of characteristic $0$.
\begin{definition}
A Cartier divisor $D$ (or invertible sheaf $\co(D)$) on $X$ is said to be numerically effective, or nef, if $D.C\geq 0$ for any closed integral curve $C$ on $X$.
\end{definition}
The next lemma summarizes some basic properties of nef divisors (see \cite{Kl}).
\begin{lemma}
Nef divisors satisfy the following:\\\\
(a).  Let $n=\dim X$.  If $D_1,\ldots,D_n$ are nef divisors on $X$ then $D_1.D_2.\ldots.D_n\geq 0$.\\
(b).  Let $D$ be a nef divisor and $A$ an ample divisor on $X$.  Then $A+D$ is ample.\\
(c).  Let $f:X \to Y$ be a morphism and let $D$ be a nef divisor on $Y$.  Then $f^*\co(D)$ is nef on $X$.  
\end{lemma}

Recall that we have defined $\kappa(D)$ and quasi-ampleness for a Cartier divisor (Definitions \ref{defk} and \ref{defbig}).  It is always true that $\kappa(D) \leq \dim X$, so $D$ is quasi-ample (or big) if and only if it has the largest possible dimension for a divisor on $X$.  For nef divisors it is possible to give a more numerical criterion for a divisor to be quasi-ample.  It is also possible in this case to get an asymptotic formula for $h^0(nD)$.  We have the following lemma, due to Sommese, as it appears in \cite{Ka}.

\begin{lemma}
\label{nefbig}
Let $D$ be a nef divisor on a nonsingular projective variety $X$.  Let $q=\dim X$.  Then $h^0(nD)=\frac{D^q}{q!}n^q+O(n^{q-1})$.  In particular, $D^q>0$ if and only if $D$ is quasi-ample.
\end{lemma}
\begin{proof}
Let $K_X$ denote the canonical divisor on $X$.  Let $L$ be an ample divisor on $X$ such that $L+K_X$ is very ample.  Since $D$ is nef, $nD+L$ is ample, and so by Kodaira's vanishing theorem we have
\begin{equation*}
H^i(X,\co(nD+L+K_X))=0 \text{ for } i>0.
\end{equation*}
Therefore,
\begin{equation*}
h^0(nD+L+K_X)=\chi(\co(nD+L+K_X))=\frac{D^q}{q!}n^q+O(n^{q-1})
\end{equation*}
by Riemann-Roch.  Let $Y$ be a general member of the linear system $|L+K_X|$, so that $Y$ is nonsingular and irreducible.  Then we have an exact sequence
\begin{equation*}
0 \to H^0(X,\co(nD)) \to H^0(X,\co(nD+L+K_X))\to H^0(Y,i^*\co(nD+L+K_X))
\end{equation*}
where $i:Y\to X$ is the inclusion map.  But since $\dim Y=q-1$, we have $\dim H^0(Y,i^*\co(nD+L+K_X))\leq O(n^{q-1})$.  It follows that $h^0(nD)=\frac{D^q}{q!}n^q+O(n^{q-1})$.
\end{proof}
Since we will use it multiple times, we state the exact sequence used above as a lemma.

\begin{lemma}
\label{exact}
Let $D$ be an effective Cartier divisor on $X$ with inclusion map $i:D \to X$.  Let $E$ be any Cartier divisor on $X$.  Then we have exact sequences
\begin{align}
&0 \to \co(E-D) \to \co(E) \to i_{*}(i^*\co(E)) \to 0\\
&0 \to H^0(X,\co(E-D))\to H^0(X,\co(E)) \to H^0(D,i^*(\co(E)).
\end{align}
\end{lemma}
\begin{proof}
If $D$ is an effective Cartier divisor, then a fundamental exact sequence is
\begin{equation*}
0 \to \co(-D) \to \mathcal{O}_X \to i_{*} \mathcal{O}_D \to 0.
\end{equation*}
Tensoring with $\co(E)$ and using the projection formula, we get the first exact sequence.  Taking global sections then gives the second exact sequence.
\end{proof}
We can prove a little more for surfaces.
\begin{lemma}
\label{surfbig}
Let $D$ be an effective divisor on a nonsingular projective surface $X$.  If $D^2>0$ then $h^0(nD)\geq \frac{n^2D^2}{2}+O(n)$ and $D$ is quasi-ample.
\end{lemma}
\begin{proof}
By Riemann-Roch,
\begin{equation}
h^0(nD)-h^1(nD)+h^0(K-nD)=\frac{n^2D^2}{2}-\frac{nD.K}{2}+1+p_a.
\end{equation}
Since $D$ is effective, $D \neq 0$, $h^0(K-nD)=0$ for $n\gg 0$ (for example, choose $n>K.H$ where $H$ is an ample divisor).  We also have $h^1(nD)\geq 0$, so $h^0(nD)\geq \frac{n^2D^2}{2}+O(n)$ and $D$ is quasi-ample. 
\end{proof}
It is not always true that if $D$ is nef then $h^0(E-D) \leq h^0(E)$.  If $h^0(D)=0$  (for example if $D$ corresponds to a non-zero torsion element of Pic $X$) then when $E=D$ we have $h^0(E-D)=h^0(D-D)=h^0(0)=1 > h^0(E)=0$.  We will want some control over $h^0(E-D)$ when $D$ is nef, and so we prove the following weak lemma.

\begin{lemma}
\label{nef}
Let $X$ be a nonsingular projective variety of dimension $q$.  Let $D$ be a nef divisor on $X$.  Let $E$ be any divisor on $X$.  Then
\begin{equation*}
h^0(nE-mD) \leq h^0(nE)+O(n^{q-1})
\end{equation*}
for all $m,n\geq 0$, where the implied constant is independent of $m$.
\end{lemma}
\begin{proof}
We first claim that if $F$ is any nef divisor then there exists a divisor $C$, independent of $F$, such that $h^0(C+F)>0$.  Explicitly, we may take $C=(q+2)A+K_X$, where $A$ is a very ample divisor on $X$.  We prove this by induction on the dimension $q$.  The case $q=1$ is easy.  For the inductive step, we have an exact sequence
\begin{multline*}
0\to H^0(X,\co((q+1)A+K_X+F))\to H^0(X,\co((q+2)A+K_X+F)) \\
\to H^0(Y,i^*(\co((q+2)A+K_X+F)))\to H^1(X,\co((q+1)A+K_X+F))
\end{multline*}
where $Y$ is an irreducible nonsingular element of $|A|$ with inclusion map $i:Y\to X$.  Since $(q+1)A+F$ is ample, by Kodaira vanishing, the last term above is $0$.  Since $\omega_Y\cong i^*(\co(A+K_X))$, by induction we get that $\dim H^0(Y,i^*(\co((q+2)A+K_X+F)))>0$.  Since the penultimate map in the exact sequence above is surjective, we therefore also have $h^0((q+2)A+k_X+F)=h^0(C+F)>0$, proving our claim.
Then we have
\begin{equation*}
h^0(nE-mD)\leq h^0(nE-mD+(C+mD))= h^0(nE+C) \leq h^0(nE) +O(n^{q-1})
\end{equation*}
independently of $m$, where the last inequality follows from Lemma \ref{exact} as in the proof of Lemma \ref{nefbig}.
\end{proof}

\section{Fundamental Theorems on Large Divisors}
\label{smain}
In this section we prove a slightly expanded version of a theorem of Corvaja and Zannier and its analogue for holomorphic curves.  These theorems will be fundamental to our future results.

Let $D$ be a divisor on a nonsingular projective variety $X$ defined over a field $k$.  Let $\kk(X)$ denote the function field of $X$ over $\kk$.  We will write $D\geq E$ if $D-E$ is effective.  Let div$(f)$ denote the principal divisor associated to $f$.  Let $L(D)$ be the $\kk$-vector space $L(D)=\{f \in \kk(X)|\text{div}(f)\geq -D\}$, and let $l(D)=\dim L(D)=h^0(D)$.  If $E$ is a prime divisor we let $\text{ord}_E f$ denote the coefficient of $E$ in div$(f)$.  We make the following definition.
\begin{definition}
Let $D$ be an effective divisor on a nonsingular projective variety $X$ defined over a field $k$.  Then we define $D$ to be a very large divisor on $X$ if for every $P\in D(\kk)$ there exists a basis $B$ of $L(D)$ such that $\text{ord}_E\prod_{f \in B}f>0$ for every irreducible component $E$ of $D$ such that $P\in E$.  We define $D$ to be a large divisor if some nonnegative integral linear combination of its irreducible components is very large on $X$. 
\end{definition}
\begin{remark}
\label{remlarge}
Suppose $D$ is very large.  Let $P\in D$ and let $\mathcal E$ be the set of irreducible components $E$ of $D$ such that $P\in E$.  If $B$ is a basis of $L(D)$ that has the property in the definition of very large with respect to $P$, then $B$ also works as a basis with respect to any $Q\in \bigcap_{E \in \mathcal{E}}E$.  Thus, it is easily seen that in the definition of very large one only needs to use bases $B\in \mathcal{B}$ for some finite set of bases $\mathcal{B}$ for any very large divisor $D$.
\end{remark}
We will see (Theorem \ref{cor3}) for example that on any nonsingular projective variety $X$ the sum of sufficiently many ample effective divisors in general position is large.  On the other hand, it is obvious from the definition that if $D$ is an irreducible effective divisor on $X$ then $D$ cannot be large.  Roughly speaking, large divisors have a lot of irreducible components of high $D$-dimension.
With this definition we have the following theorems.
\renewcommand{\thetheorem}{\arabic{section}.\arabic{theorem}A}
\setcounter{theoremb}{\value{theorem}}
\begin{theorema}[Corvaja-Zannier]
\label{maina}
Let $X$ be a nonsingular projective variety defined over a number field $k$.  Let $S\subset M_k$ be a finite set of places of $k$ containing the archimedean places.  Let $D$ be a large divisor on $X$ defined over $k$.  Then there does not exist a Zariski-dense set of $D$-integral points on $X$.  Furthermore, if $D$ is very large and $\Phi_D$ is a rational map to projective space corresponding to $D$, then there exists a proper closed subset $Z\subset X$ depending only on $D$ (and not $k$ or $S$) such that $\Phi_D(R\backslash Z)$ is finite for any set $R$ of $D$-integral points on $X$.  In particular, if $\Phi_D$ is birational, $X\backslash D$ is quasi-Mordellic.
\end{theorema}
\begin{theoremb}
\label{mainb}
Let $X$ be a nonsingular complex projective variety.  Let $D$ be a large divisor on $X$.  Then there does not exist a holomorphic map $f:\mathbb{C} \to X \backslash D$ with Zariski-dense image.  Furthermore, if $D$ is very large and $\Phi_D$ is a rational map to projective space corresponding to $D$, then there exists a proper closed subset $Z\subset X$ depending only on $D$ such that for all holomorphic maps  $f:\mathbb{C} \to X \backslash D$, either $f(\mathbb{C})\subset Z$ or $\Phi_D\circ f$ is constant.  In particular, if $\Phi_D$ is birational, $X\backslash D$ is quasi-Brody hyperbolic.
\end{theoremb}
Theorem \ref{maina} appears, essentially, in the proof of the Main Theorem in \cite{Co2}, and for curves in \cite{Co}.  We have added the last two statements to the theorem by using Vojta's result on the exceptional hyperplanes in the Schmidt Subspace Theorem.

Given these theorems, many of our results mentioned in the introduction reduce to showing that certain divisors are large.
Let us prove Theorem \ref{maina} first.  Before proving this theorem, we need a lemma.
\renewcommand{\thetheorem}{\arabic{section}.\arabic{theorem}}
\begin{lemma}
\label{seq}
Let $X$ be a projective variety defined over a number field $k$.  Let $R\subset X(k)$ be a Zariski-dense subset of $X$.  Let $v\in M_k$.  Then there exists a point $P$ in $X(k_v)$ and a sequence $\{P_i\}$ in $R$ such that $\{P_i\}\to P$ in the $v$-topology on $X(k_v)$ and $\bigcup \{P_i\}$ is Zariski-dense in $X$.
\end{lemma}
\begin{proof}
We will always be working in the $v$-topology on $X(k_v)$.  First we claim that there exists a $P$ in $\overline{R}\subset X(k_v)$ such that for every neighborhood $U$ of $P$ in $X(k_v)$, $U\cap R$ is Zariski-dense in $X$.  Indeed, suppose there is no such $P$.  Then for each $P$ in $R$, let $U_P$ be a neighborhood of $P$ such that $U_P \cap R$ is not Zariski-dense in $X$.  Since $X(k_v)$ is compact because $X$ is projective, $\overline{R}$ is compact, so we may cover $\overline{R}$ by finitely many open sets $U_{P_1},\ldots,U_{P_n}$.  But then $R=(U_{P_1}\cap R)\cup\cdots \cup (U_{P_n}\cap R)$ is not Zariski-dense in $X$, a contradiction.

Now pick some $P$ as in the claim above.  Embed $X$ in $\mathbb{P}^n_k$ for some $n$.  Since $k$ is countable, the set of hypersurfaces in $\mathbb{P}^n_k$ not containing $X$ is countable.  Let $\{H_i\}$ be an enumeration of these.  There also exists a countable collection of neighborhoods $\{U_i\}$ of $P$ in $X(k_v)$ such that $U_i \subset U_j$ for $i>j$ and $\bigcap U_i=\{P\}$.  Since $U_i\cap R$ is Zariski-dense in $X$, for all $i$ there exists a $P_i \in U_i \cap R$ such that $P_i \notin H_i$.  Then $\{P_i\}\to P$ in $X(k_v)$ and $\bigcup \{P_i\}$ is Zariski-dense in $X$ since it is not contained in any hypersurface.
\end{proof}
\begin{proof}[Proof of Theorem \ref{maina}]
Let $D$ be a large divisor and $S$ and $X$ as in Theorem \ref{maina}.  Clearly, we may reduce to the case where $D$ is very large.  Extending $k$ if necessary and enlarging $S$, we may assume without loss of generality that every irreducible component of $D$ is defined over $k$ and that all of the finitely many functions in $L(D)$ we use (see Remark \ref{remlarge}) are defined over $k$.  Let $\{\phi_1,\ldots,\phi_{l(D)}\}$ be a basis of $L(D)$ over $k$.  Let $R$ be a $(D,S)$-integral set of points on $X$.  It suffices to prove the theorem in the case that $\overline{R}$ is irreducible.  By repeatedly applying Lemma~\ref{seq}, we see that there exists a sequence $P_i$ in $R$ such that for each $v$ in $S$, $\{P_i\}$ converges to a point $P_v\in X(k_v)$ and $\bigcup \{P_i\}$ is Zariski-dense in $\overline{R}$.

Let $S'$ be the set of places $v\in S$ such that $P_v\in D(k_v)$, and let $S''=S\backslash S'$.  Since $D$ is very large, for each $v\in S'$ let $L_{iv},i=1,\ldots,l(D)$ be a basis for $L(D)$ such that $\text{ord}_E\prod_{i=1}^{l(D)}L_{iv}>0$ for all irreducible components $E$ of $D$ such that $P_v\in E(k_v)$.  Of course each $L_{iv}$ is a linear form in the $\phi_j$'s over $k$.
For $v\in S''$, we set $L_{jv}=\phi_j$ for $j=1,\ldots,l(D)$.  Let $\phi(P)=(\phi_1(P),\ldots,\phi_{l(D)}(P))$ for $P\in X\backslash D$.  Let $H_{jv}$ denote the hyperplane in $\mathbb{P}^{l(D)-1}$ determined by $L_{jv}$ with respect to the basis $\phi_1,\ldots,\phi_{l(D)}$.  Let $\lambda_{H_{jv},v}$ be the Weil function for $H_{jv}$ given in Equation (\ref{Weila}).  
We will now show that there exists $\epsilon>0$ and a constant $C$ such that
\begin{equation}
\label{Schmidt}
\sum_{v\in S}\sum_{j =1}^{l(D)} \lambda_{H_{jv},v}(\phi(P_i)) > (l(D)+\epsilon)h(\phi(P_i))+C.
\end{equation}
Since $R$ is a set of $(D,S)$-integral points, we have 
\begin{equation*}
h(\phi(P_i))<\sum_{v \in S}\log{\max_j|\phi_j(P_i)|_v}+O(1).
\end{equation*}
Using this it suffices to prove that
\begin{equation*}
\sum_{v\in S}\sum_{j =1}^{l(D)}\log \max_{j'} \frac{|\phi_{j'}(P_i)|_v}{|L_{jv}(P_i)|_v}>(l(D)+\epsilon)\sum_{v\in S}\log{\max_{j'}|\phi_{j'}(P_i)|_v}+C'
\end{equation*}
for some $C'$ or rearranging things, simplifying, and exponentiating
\begin{equation*}
\prod_{v\in S} |\max_{j'}(\phi_{j'}(P_i))^{\epsilon} \prod_{j=1}^{l(D)}L_{jv}(P_i)|_v
\end{equation*}
is bounded for some $\epsilon>0$.  Let 
\begin{equation*}
M=\max\{-\text{ord}_E \phi_j|E \text{ is an irreducible component of } D, j=1,\ldots,l(D)\}.
\end{equation*}
Let $\epsilon=\frac{1}{M}$.  For $v\in S''$ both $|\phi_{j'}(P_i)|_v$ and $|L_{jv}(P_i)|_v$ are bounded for all $i$ since $P_v \notin D(k_v)$ and $\phi_{j'}$ and $L_{jv}$ have poles lying only in the support of $D$.  Let $v\in S'$.  So $P_v\in D(k_v)$.  It follows from the definition of $M$ and the fact that $\text{ord}_E \prod_{i=1}^{l(D)}L_{iv}>0$ for any irreducible component $E$ of $D$ such that $P_v \in E(k_v)$ that $\text{ord}_E\phi_{j'} (\prod_{i=1}^{l(D)}L_{iv})^M\geq -M+ M \geq 0$ for any irreducible component $E$ of $D$ such that $P_v \in E(k_v)$.  Since the $\phi_{j'}$ and $L_{iv}$ have poles only in the support of $D$, it follows from the previous order computation that $|\max_{j'}(\phi_{j'}(P_i))^{\epsilon} \prod_{j=1}^{l(D)}L_{jv}(P_i)|_v$ is bounded for all $i$ and all $v\in S$ when $\epsilon=\frac{1}{M}>0$.  So we have proved Equation (\ref{Schmidt}). 

Note that either $h(\phi(P_i))\to \infty$ as $i\to \infty$ or $\phi(P_i)=\phi(\overline{R})$, and $\phi(P_i)$ is constant for all $i$.  In the latter case the theorem is proved, so we may assume the former.  Therefore, making $\epsilon$ smaller, we see that Equation (\ref{Schmidt}) holds with $C=0$ for all but finitely many $i$.  So by Schmidt's Subspace Theorem, there exists a finite union of hyperplanes $Z\subset \mathbb{P}^{l(D)-1}$ such that all but finitely many of the points in the set $\{\phi(P_i)=(\phi_1(P_i),\ldots,\phi_{l(D)}(P_i))|i\in \mathbb{N}\}$ lie in $Z$.  Using Remark \ref{remlarge} we see that we may choose the hyperplanes $H_{iv}$ used above from a finite set of hyperplanes independent of $R$.  Therefore, using the statement on the exceptional hyperplanes in the Schmidt Subspace Theorem, we see that $Z$ may be chosen to depend only on $D$ and not $R$, $k$, or $S$.   Since it was assumed that $\overline{R}$ is irreducible and $\phi(\overline{R})$ is not a point, it follows that $\phi(R)\subset Z$.  Since $\phi_1,\ldots,\phi_{d}$ are linearly independent functions in $K(X)$ and $Z$ is a finite union of hyperplanes, it follows that $\phi^{-1}(Z)$ is a finite union of proper closed subvarieties of $X$.  So $R\subset \phi^{-1}(Z)$ and the theorem is proved.
\end{proof}
The proof of Theorem \ref{mainb} is very similar.
\begin{proof}[Proof of Theorem \ref{mainb}]
Since our assertion depends only on the support of $D$ we may assume without loss of generality that $D$ is very large on $X$.  Let $f:\mathbb{C} \to X \backslash D$ be a holomorphic map.  By Remark \ref{remlarge} there exists a finite set $J$ of elements in $L(D)$ such that for any $P\in D$ there exists a subset $I \subset J$ that is a basis of $L(D)$ such that $\text{ord}_E \prod_{g\in I}g>0$ for every irreducible component $E$ of $D$ such that $P\in E$.  Let $\phi_1,\ldots,\phi_{l(D)}$ be a basis for $L(D)$.  Let $\phi=(\phi_1,\ldots,\phi_{l(D)}):X\backslash D \to \mathbb{P}^{l(D)-1}$.  Let $J'$ be the set of linear forms $L$ in $l(D)$ variables over $\mathbb{C}$ such that $L\circ \phi \in J$.  If $L$ is a linear form, let $H_L$ be the corresponding hyperplane.  We will now show that there exists $\epsilon>0$ and a constant $C$ such that
\begin{equation}
\label{Cartan}
\int_{0}^{2\pi} \max_I \sum_{L \in I} \lambda_{H_L}(\phi \circ f(re^{i\theta}))\frac{d\theta}{2\pi} > (l(D)+\epsilon)T_{\phi \circ f}(r)-C
\end{equation}
for all $r>0$, where the max is taken over subsets $I \subset J'$ such that $I$ consists of exactly $l(D)$ linearly independent linear forms.  Substituting the definition of the Weil function in Equation (\ref{Weilb}) and the definition of $T_{\phi \circ f}$, after some manipulation the inequality in Equation (\ref{Cartan}) becomes
\begin{equation*}
\int_{0}^{2\pi} \epsilon \log|\phi\circ f(re^{i\theta})|+\min_I \sum_{L \in I} \log |L\circ \phi \circ f(re^{i\theta})| \frac{d\theta}{2\pi}<C
\end{equation*}
with $I$ as before.  Since $|\phi\circ f(re^{i\theta})|\leq \sqrt{l(D)} \max_j |\phi_j \circ f(re^{i\theta})|$ it clearly suffices to show that
\begin{equation}
\label{Cartan2}
\max_j|\phi_j \circ f(re^{i\theta})|^{\epsilon} \min_I \prod_{L \in I} |L\circ \phi \circ f(re^{i\theta})|
\end{equation}
is bounded independently of $r$ and $\theta$ for some $\epsilon>0$.
Let $D_1,\ldots,D_m$ be the irreducible components of $D$.  Let
\begin{equation*}
M=\max\{-\text{ord}_{D_i} \phi_j|i=1,\ldots,m, j=1,\ldots,l(D)\}.
\end{equation*}
We will work in the classical topology.  Let $P\in D$.  Then there exists a neighborhood $U$ of $P$ such that for all $Q\in \overline{U}$ if $Q \in D_i$ for some $i$ then $P \in D_i$.  Let $I\subset J'$ be a subset of $J'$ such that $\text{ord}_{D_i} \prod_{L\in I} L \circ \phi>0$ for all $i$ such that $P\in D_i$.  If $P \in D_i$, then by the definition of $M$ we have $\text{ord}_{D_i}\phi_j (\prod_{L\in I} L \circ \phi)^M\geq 0$ for all $j$.  By the definition of $U$ we see that $\phi_j (\prod_{L\in I} L \circ \phi)^M$ is bounded for all $j$ on the compact set $\overline{U}$.  Since $D$ is compact and may be covered by such sets we see that $\max_j|\phi_j|\min_I \prod_{L \in I} |L\circ \phi|^M$ is bounded on $X\backslash D$ (using also that away from $D$ everything is obviously bounded since the $\phi_j$'s have poles only in $D$).  Therefore Equation (\ref{Cartan2}) is bounded independently of $r$ and $\theta$ for $\epsilon=\frac{1}{M}$.  

If $\phi \circ f$ is constant then there is nothing to prove, so assume otherwise.  Then  $T_{\phi \circ f}(r)\to \infty$ as $r\to \infty$, and so making $\epsilon$ smaller, we see that we have proven the inequality (\ref{Cartan}) with $C=0$ for all sufficiently large $r$.  Therefore by Cartan's Second Main Theorem, there exists a finite union of hyperplanes $Z\subset \mathbb{P}^{l(D)-1}$ depending only on $D$ (the $H_L$'s depended only on $D$) such that $\phi(f(\mathbb{C}))\subset Z$.  Since the $\phi_j$'s are linearly independent and $Z$ is a finite union of hyperplanes, $\phi^{-1}(Z)$ is a finite union of closed subvarieties of $X$ and $f(\mathbb{C})\subset \phi^{-1}(Z)$.  
\end{proof}
\begin{remark}
If $D$ is very large and one can explicitly compute the map $\phi$ and the hyperplanes used in the above proofs, then one can explicitly compute the closed set $Z$ in the theorems above.  This follows from the explicit description of the exceptional hyperplanes in \cite{Vo6} and \cite{Vo3}.
\end{remark}

\section{Large Divisors}
For an effective divisor $D=\sum_{i=1}^r D_i$ on $X$ and $P\in D(\kk)$, we let $D_P=\sum_{i:P\in D_i}D_i$.
\begin{lemma}
\label{large}
Let $D=\sum_{i=1}^r D_i$ be a divisor on a nonsingular projective variety $X$ with $D_i$ effective for each $i$.  Let $P\in D$.  Let $f_P(m,n)=l(nD-mD_P)-l(nD-(m+1)D_P)$.  If there exists $n>0$ such that $\sum_{m=0}^{\infty}(m-n)f_P(m,n)>0$ for all $P\in D$ then $nD$ is very large.
\end{lemma}
\begin{proof}
Let $n>0$ be such that $\sum_{m=0}^{\infty}(m-n)f_P(m,n)>0$ for all $P\in D$.  This sum is clearly finite for all $P\in D$ and we let $M_P(n)$ be the largest integer such that $f_P(M_P(n),n)>0$.  Let $P\in D$.  Let $M=M_P(n)$.  Let $V_j=L(nD- jD_P)$.  So $\dim V_j/V_{j+1}=f_P(j,n)$.  We have $L(nD)=V_0 \supset V_1 \supset \ldots \supset V_M\neq 0$.  Choose a basis of $V_M$ and successively complete it to bases of $V_{M-1},V_{M-2},\ldots,V_0$, to obtain a basis $f_1,\ldots,f_{l(nD)}$.  Let $E$ be an irreducible component of $D$ such that $P \in E$. If $f_j \in V_m$ then $\text{ord}_E f_j\geq (m-n)\ord_ED$.  So we get that $\text{ord}_E\prod_{i=1}^{l(nD)}f_i\geq (\ord_ED)\sum_{m=0}^{M}(m-n)f_P(m,n)>0$.  So $nD$ is very large.
\end{proof}
\begin{theorem}
\label{cor2}
Let $X$ be a nonsingular projective variety.  Let $q= \dim X$.  Let $D=\sum_{i=1}^{r}D_i$ be a divisor on $X$ such that $D_i$ is effective and nef for each $i$.  Suppose also that every irreducible component of $D$ is nonsingular.  If
\begin{equation*}
D^q>2q D^{q-1}.D_P, \qquad \forall P\in D
\end{equation*}
then $nD$ is very large for $n\gg 0$.  In particular, $D$ is large.
\end{theorem}
\begin{proof}
Let $P \in D$.  Let $D_P=\sum_{j=1}^{k}a_j E_j$, where each $E_j$ is a distinct prime divisor.  Repeatedly applying Lemma~\ref{exact}, we obtain
\begin{multline*}
\dim H^0(X,\co(nD-m D_P))-\dim H^0(X,\co(nD-(m+1)D_P)) \\
\leq \sum_{j=1}^k \sum_{l=0}^{a_{j}-1}\dim H^0(E_{j},i^*_{E_{j}}\co(nD-mD_P -\sum_{j'=1}^{j-1}a_{j'}E_{j'}-lE_{j}))
\end{multline*}
It follows from the fact that $D_P$ is nef, Lemma \ref{exact}, and Lemma \ref{nef} that
\begin{multline*}
\dim H^0(E_{j},i^*_{E_{j}}\co(nD-mD_P -\sum_{j'=1}^{j-1}a_{j'}E_{j'}-lE_{j}))\\
\leq \dim H^0(E_{j},i^*_{E_{j}}\co(nD))+O(n^{q-2}).
\end{multline*}
Therefore,
\begin{multline*}
\dim H^0(X,\co(nD-m D_P))-\dim H^0(X,\co(nD-(m+1)D_P))\\
\leq \sum_{j=1}^k a_j \dim H^0(E_j,i_{E_j}^*\co(nD))+O(n^{q-2}).
\end{multline*}
Since $D$ is nef, $l(nD)=\frac{n^q}{q!}D^q + O(n^{q-1})$.  Since $i_{E_j}^*\co(D)$ is also nef, we have $\dim  H^0(E_j,i_{E_j}^* \co(nD))= \frac{n^{q-1}}{(q-1)!} D^{q-1}.E_j+ O(n^{q-2})$.  So 
\begin{equation*}
f_P(m,n)\leq \frac{n^{q-1}}{(q-1)!} \sum_{j=1}^k a_j D^{q-1}.E_j+ O(n^{q-2}) =\frac{n^{q-1}}{(q-1)!} D^{q-1}.D_P+ O(n^{q-2}).
\end{equation*}
To use this estimate, we borrow a lemma from \cite{Co2}.
\begin{lemma}
\label{CZlemma}
Let $h$ and $R$ be integers with $R\leq h$ and let $x_1,\ldots,x_h,U_1,\ldots,U_R$ be real numbers.  If $0\leq x_i\leq U_i$ for $i=1,\ldots, R$ and $\sum_{j=1}^RU_j\leq \sum_{j=1}^hx_j$ then $\sum_{j=1}^hjx_j\geq \sum_{j=1}^RjU_j$.
\end{lemma}
\begin{proof}
We have
\begin{align*}
\sum_{j=1}^RjU_j+\sum_{j=1}^h(R+1-j)x_j&\leq \sum_{j=1}^RjU_j+\sum_{j=1}^R(R+1-j)x_j\\
&\leq \sum_{j=1}^RjU_j+\sum_{j=1}^R(R+1-j)U_j=(R+1)\sum_{j=1}^RU_j
\end{align*}
So, rearranging things.
\begin{equation*}
\sum_{j=1}^h jx_j\geq \sum_{j=1}^RjU_j+(R+1)\left(\sum_{j=1}^hx_j-\sum_{j=1}^RU_j\right)
\end{equation*}
and the last term is positive by assumption.
\end{proof}
Let $R_n=\frac{n^q}{q!}D^q$ and $S_n= \frac{n^{q-1}}{(q-1)!} D^{q-1}.D_P$.  In the notation of Lemma \ref{large}, we have
\begin{equation*}
\sum_{m=0}^{M_P(n)}f_P(m,n)=l(nD)=R_n+O(n^{q-1}).
\end{equation*}
and $f_P(m,n)\leq S_n+O(n^{q-2})$.  We will assume from now on that $S_n\neq 0$ (the case $S_n=0$ is similar).
Then using our estimate, we have $M_P(n) \geq \frac{R_n}{S_n}+O(1)$ and $\sum_{m=0}^{\frac{R_n}{S_n}+O(1)}(S_n+ O(n^{q-2}))\leq \sum_{m=0}^{M_P(n)}f_P(m,n)$.  So using Lemma \ref{CZlemma}, for $n \gg 0$ we get the estimate
\begin{align*}
\sum_{m=0}^{M_P(n)}(m-n)f_P(m,n) & \geq \sum_{m=0}^{\frac{R_n}{S_n}+O(1)}m(S_n+O(n^{q-2}))-n\sum_{m=0}^{M_P(n)}f_P(m,n)\\
&\geq \frac{R_n^2}{2S_n}-nR_n+O(n^q)\\
&\geq \frac{R_n}{S_n}\left[\frac{n^q}{2q!}\left(D^q-2q D^{q-1}.D_P\right)+O(n^{q-1})\right].
\end{align*}
So for $n \gg 0$, $ \sum_{m=0}^{M_P}(m-n)f_P(m,n)>0$ if $D^q>2q D^{q-1}.D_P$.  Then we are done by Lemma \ref{large}.
\end{proof}
When $q=1$ we obtain
\begin{corollary}
Let $D$ be an effective divisor on a nonsingular projective curve $X$.  If $D$ is a sum of more than 2 distinct points on $X$ then $D$ is large. 
\end{corollary}
By Theorems \ref{maina} and \ref{mainb} we then recover
\renewcommand{\thetheorem}{\arabic{section}.\arabic{theorem}A}
\setcounter{theoremb}{\value{theorem}}
\begin{corollarya}
Siegel's theorem (Theorem~\ref{Siegel2})
\end{corollarya}
\begin{corollaryb}
Picard's theorem (Theorem \ref{Picard})
\end{corollaryb}
\renewcommand{\thetheorem}{\arabic{section}.\arabic{theorem}}
Actually we have only proved these theorems for nonsingular curves $\tilde{C}$.  However, the general case follows from this case by looking at the normalization of $\tilde{C}$.

Suppose that we have a divisor $D=\sum_{i=1}^{r}D_i$ satisfying the hypotheses of Theorem~\ref{cor2}.  We would like to modify $D$ to a divisor $D'=\sum_{i=1}^{r}a_iD_i$ so that we may optimally apply the theorem.  When each $D_i$ is ample, this amounts to choosing the $a_i$'s so that in the embedding given by $nD'$ for $n\gg 0$ the degree of each $a_iD_i$ is the same.  In terms of intersection theory, we would like $a_iD_i.(D')^{q-1}$ to be the same for each $i$.  We make the following definition:
\begin{definition}
Let $X$ be a nonsingular projective variety.  Let $q=\dim X$.  Let $D=\sum_{i=1}^rD_i$ be a divisor on $X$ with $D_1,\ldots,D_r$ effective..  Then $D$ is said to have equidegree with respect to $D_1,\ldots,D_r$ if $D_i.D^{q-1}=\frac{D^q}{r}$ for $i=1,\ldots,r$.  We will say that $D$ is equidegreelizable (with respect to $D_1,\ldots,D_r$) if there exist real numbers $a_i>0$ such that if $D'=\sum_{i=1}^ra_iD_i$ then $D'$ has equidegree with respect to $a_1D_1,\ldots,a_r D_r$.  (extending intersections to $\Div X\otimes \mathbb{R}$ in the canonical way).
\end{definition}
We will frequently omit the reference to the $D_i$'s when it is clear what we mean.
\begin{lemma}
\label{equi}
Let $X$ be a nonsingular projective variety.  Let $q=\dim X$.  Let $D_1,\ldots,D_r$ be divisors on $X$ with $D_i^q>0$ for all $i$.  Suppose that all $q$-fold intersections of the $D_i$'s are nonnegative.  Then $\sum_{i=1}^r D_i$ is equidegreelizable with respect to $D_1,\ldots,D_r$.
\end{lemma}
\begin{proof}
Consider the function $f(a_1,\ldots,a_r)=(\sum_{i=1}^r e^{a_i}D_i)^q$ subject to the constraint $g(a_1,\ldots,a_r)=\sum_{i=1}^r a_i=0$.  Since all $q$-fold intersections of the $D_i$'s are nonnegative, $f(a_1,\ldots,a_r)\geq e^{q a_i}D_i^q$ for any $i$.  Since $D_i^q>0$ for all $i$, as $\max \{a_i\}\to \infty$ we have $f(a_1,\ldots,a_r) \to \infty$.  It follows that $f$ attains a minimum on the plane $\sum_{i=1}^r a_i=0$.  Therefore there exists a solution $\lambda,a_1,\ldots,a_r$ to the Lagrange multiplier equations $g=0,\frac{\partial f}{\partial a_i}=e^{a_i}D_i.(\sum_{i=1}^r e^{a_i} D_i)^{q-1}=\lambda \frac{\partial g}{\partial a_i}=\lambda, i=1,\ldots,r$.  So $D'=\sum_{i=1}^r e^{a_i}D_i$ has equidegree with respect to $D_1,\ldots,D_r$ and trivially $e^{a_i}>0$.  
\end{proof}
We give an example to show that not all divisor sums are equidegreelizable.
\begin{example}
Let $X=\mathbb{P}^1 \times \mathbb{P}^1$.  Let $D_1=P_1 \times \mathbb{P}^1,D_2=P_2 \times \mathbb{P}^1$, and $D_3=\mathbb{P}^1 \times Q$, where $P_1,P_2$, and $Q$ are points in the various $\mathbb{P}^1$'s.  So $D_1.D_2=D_1^2=D_2^2=D_3^2=0$ and $D_1.D_3=D_2.D_3=1$.  Let $D=a_1D_1+a_2D_2+a_3D_3$.  Since $a_3D_3.D=a_1D_1.D+a_2D_2.D$, it is clear that there do not exist $a_1,a_2,a_3>0$ such that $a_iD_i.D=\frac{D^2}{3}$ for $i=1,2,3$.  So $D=D_1+D_2+D_3$ is not equidegreelizable with respect to $D_1, D_2$, and $D_3$.
\end{example}
With the above definition, we have the following theorem.
\begin{theorem}
\label{cor3}
Let $X$ be a nonsingular projective variety.  Let $q=\dim X$.  Let $D=\sum_{i=1}^r D_i$ be a quasi-ample divisor on $X$ equidegreelizable with respect to $D_1,\ldots, D_r$, with $D_1,\ldots, D_r$ nef and effective.  Suppose that every irreducible component of $D$ is nonsingular.  Suppose that the intersection of any $m+1$ distinct $D_i$'s is empty.  If $r>2mq$ then $D$ is large.  Furthermore, there exists a very large divisor $E$ with the same support as $D$ such that $\Phi_E$ is birational.
\end{theorem}
\begin{proof}
Since $D$ is equidegreelizable, we may find positive integers $a_i$ such that if $D'=\sum_{i=1}^{r}a_iD_i$ then $\frac{a_iD_i.(D')^{q-1}}{(D')^q}$ is arbitrarily close to $\frac{1}{r}$ for each $i$.  Note that $D'$ is again quasi-ample.  Since for any $P\in D(\kk)$, $P$ belongs to at most $m$ divisors $D_i$, and $r>2mq$, we have that
\begin{equation*}
2q(D')^{q-1}.(D')_P=2q \sum_{i: P\in D_i(\kk)}a_i D_i.(D')^{q-1}<(D')^q.
\end{equation*}
So the hypotheses of Theorem~\ref{cor2} are satisfied and so $nD'$ is very large for $n \gg 0$.  The last statement then follows from the fact that $D'$ is quasi-ample.
\end{proof}
\begin{lemma}
\label{reduce}
Let $X$ be a complex projective variety.  Let $D=\sum_{i=1}^rD_i$ be a sum of effective Cartier divisors on $X$.  Then there exists a nonsingular projective variety $X'$, a birational morphism $\pi:X'\to X$, and a divisor $D'=\sum_{i=1}^rD_i'$ on $X'$ such that $\supp D_i'\subset \supp \pi^*D_i$ for all $i$, every irreducible component of $D'$ is nonsingular, $|D_i'|$ is base-point free for all $i$ (in particular $D_i'$ is nef), and $\kappa(D_i')=\kappa(D_i)=\dim \Phi_{D_i'}(X')$ for all $i$.
\end{lemma}
\begin{proof}
Taking a resolution of the singularities of $X$ and of the embedded singularities of the irreducible components of $D$ we may assume that $X$ and every irreducible component of $D$ are nonsingular.  For each $i$, let $m_i>0$ be such that $\dim \Phi_{m_iD_i}(X)=\kappa(D_i)$.  Let $\pi:X'\to X$ be the map obtained by blowing up the base-points of all the linear systems $|m_iD_i|$.  Then $\pi^*(m_iD_i)=D_i'+F_i$ for each $i$, where $|D_i'|$ is base-point free and $F_i$ is the fixed part of $|\pi^*(m_iD_i)|$.  We have, trivially from the definition, $\kappa(D_i)=\kappa(m_iD_i)$.  Further,  $\kappa(m_iD_i)=\kappa(\pi^*(m_iD_i))$ (in fact $l(mD_i)=l(\pi^*(mD_i))$ for all $m$ follows easily from $\pi_*\co_{X'}=\co_X$ and the projection formula). Finally, $\kappa(\pi^*(m_iD_i))=\kappa(D_i')$ since by construction $\kappa(D_i')=\max_{n>0}\dim \Phi_{nD_i'}(X')\geq \kappa(D_i)=\kappa(\pi^*(m_iD_i))$ (the other inequality being trivial).  So $\kappa(D_i')=\kappa(D_i)$ for all $i$ and therefore $X',\pi$, and $D'=\sum_{i=1}^rD_i'$ satisfy the requirements of the lemma.
\end{proof}
We now obtain one of our main results.
\renewcommand{\thetheorem}{\arabic{section}.\arabic{theorem}A}
\setcounter{theoremb}{\value{theorem}}
\begin{theorema}
\label{cor4}
Let $X$ be a projective variety defined over a number field $k$.  Let $q=\dim X$.  Let $D=\sum_{i=1}^{r} D_i$ be a divisor on $X$ defined over $k$ such that the $D_i$'s are effective Cartier divisors and the intersection of any $m+1$ distinct $D_i$'s is empty.\\\\
(a).  If $D_i$ is quasi-ample for each $i$ and $r> 2mq$ then $X\backslash D$ is quasi-Mordellic.\\
(b).  If $D_i$ is ample for each $i$ and $r> 2mq$ then $X\backslash D$ is Mordellic.
\end{theorema}
\begin{theoremb}
\label{cor4b}
Let $X$ be a complex projective variety.  Let $q=\dim X$.  Let $D=\sum_{i=1}^{r} D_i$ be a divisor on $X$ such that the $D_i$'s are effective Cartier divisors and the intersection of any $m+1$ distinct $D_i$'s is empty.\\\\
(a).  If $D_i$ is quasi-ample for each $i$ and $r>2mq$ then $X\backslash D$ is quasi-Brody hyperbolic.\\
(b).  If $D_i$ is ample for each $i$ and $r> 2mq$ then $X\backslash D$ is complete hyperbolic and hyperbolically imbedded in $X$.  In particular, $X\backslash D$ is Brody hyperbolic.
\end{theoremb}
Aside from the statement about being complete hyperbolic and hyperbolically imbedded, the same proof works for both Theorems \ref{cor4} and \ref{cor4b}.
\begin{proof}
We'll prove part (a) first.  Note that if $\pi:X'\to X$ is a birational morphism and the conclusions of part (a) of the theorems hold for $\pi^*D$ on $X'$ then they hold for $D$ on $X$.  Therefore, by Lemma \ref{reduce}, we may assume (extending $k$ in the Diophantine case if necessary) that $X$ is nonsingular, every irreducible component of $D$ is nonsingular, and $D_i$ is nef for all $i$.  The statement then follows from Lemma \ref{equi}, Theorem \ref{cor3}, and Theorems \ref{maina} and \ref{mainb}.

For part (b), we note that by (a) any set of $D$-integral points (resp. the image of any holomorphic map $f:\mathbb{C} \to X \backslash D$) is not Zariski-dense.  Let $R$ be a set of $D$-integral points (resp. the image of a holomorphic map $f:\mathbb{C} \to X \backslash D$).  Let $Y$ be an irreducible component of the Zariski-closure of $R$.  Suppose $\dim Y>0$.  Then $D$ pulls back to a sum of $r$ ample (hence quasi-ample) divisors on $Y$ such that the intersection of any $m+1$ of them is empty.  But $R\cap Y$ is a dense set of $D|_Y$-integral points on $Y$ (resp. the image of a holomorphic map $f:\mathbb{C} \to Y \backslash D$), contradicting part (a) proven above since $r>2mq>2m\dim Y$.  Therefore $\dim Y=0$.

To prove the extra hyperbolicity statements in (b) in the analytic case, we use Theorem \ref{hyperbolic}.  Let $\emptyset \subset J \subset \{1,\ldots,r\}$.  Let $s=\#J$.  Let $X'=\bigcap_{j\in J}D_j$.  We may clearly assume that $X'\not\subset D_i$ for any $i\in I\backslash J$ and that $\dim X'>0$.  Let $D'=\sum_{i\in I\backslash J}D_i|_{X'}$.  Then $D'$ is a sum of $r-s$ ample divisors on $X'$ and the intersection of any $m-s+1$ of the ample divisors is empty since $X'$ is already an intersection of $s$ of the $D_i$'s.  But $r>2mq$ implies that $r-s>2(m-s)\dim X'$.  Therefore by what we have proven above, $X'\backslash D'$ is Brody hyperbolic.  So by Theorem~\ref{hyperbolic}, $X\backslash D$ is complete hyperbolic and hyperbolically imbedded in $X$.
\end{proof}
We can prove our Main Conjectures in the simple case $m=1$ by reducing to Siegel's and Picard's theorems.  We will need the following Bertini theorem (see \cite[Th. 7.19]{Ii}).
\renewcommand{\thetheorem}{\arabic{section}.\arabic{theorem}}
\begin{theorem}
\label{Bertini}
Let $|D|$ be a base-point free linear system on a nonsingular projective variety $X$ with $\dim \Phi_{D}(X)\geq 2$.  Then every member of $|D|$ is connected and a general member of $|D|$ is nonsingular and irreducible.
\end{theorem}
\begin{lemma}
\label{m1}
Suppose $D=D_1+D_2$ is a Cartier divisor on a projective variety $X$ with $\kappa(D_1)>0,\kappa(D_2)>0$ and $D_1\cap D_2=\emptyset$.  Then $\kappa(D)=\kappa(D_1)=\kappa(D_2)=1$.
\begin{proof}
By Lemma \ref{reduce}, we may assume that $X$ is nonsingular and $|D|$ is base-point free.  If $\kappa(D)\geq 2$ then $\dim \Phi_{nD}(X)\geq 2$ for some $n>0$.  But by Theorem \ref{Bertini}, every divisor in $|nD|$ is connected, contradicting  $D_1\cap D_2=\emptyset$.
\end{proof}
\end{lemma}
\begin{theorem}
\label{tm1}
The Main Conjectures, Conjectures \ref{conjmaina},B through \ref{conj2a},B,  are true if $m=1$ (i.e. $D_i\cap D_j=\emptyset$ for all $i\neq j$).
\end{theorem}
\begin{proof}
By the above lemma, it suffices to prove the conjectures when $D=\sum_{i=1}^rD_i$ with $r>2$, and $\kappa(D)=1$.  By Lemma \ref{reduce}, we may assume that $X$ is nonsingular and $D$ is base-point free.  For $n\gg 0$, $\Phi_{nD}(X)$ is a nonsingular curve $C$ and $\Phi_{nD}$ has connected fibers.  Therefore, since $D_i\cap D_j=\emptyset$ for $i\neq j$, we have $\Phi_{nD}(X\backslash D)=C\backslash\{r\text{ points}\}$.  Since $r>2$, we are done by Siegel's and Picard's theorems.
\end{proof}
\section{A Filtration Lemma}
We'll now show how some of the results in the last section may be improved by use of a linear algebra lemma on filtrations.  The idea of using this lemma, as well as its statement and proof, are taken from the paper \cite{Co2}.  Corvaja and Zannier used it to prove a result on integral points on surfaces, and it will be essential for our results on surfaces in the next section also.
\begin{lemma}
Let $V$ be a vector space of finite dimension $d$ over a field $k$.  Let $V=W_1\supset W_2\supset \cdots \supset W_h,V=W_1^*\supset W_2^*\supset \cdots \supset W_{h^*}^*$ be two filtrations on $V$.  There exists a basis $v_1,\dots, v_d$ of $V$ which contains a basis of each $W_j$ and $W_j^*$.
\end{lemma}
\begin{proof}
The proof will be by induction on $d$.  The case $d=1$ is trivial.  By refining the first filtration, we may assume without loss of generality that $W_2$ is a hyperplane in $V$.  Let $W_i'=W_i^*\cap W_2$.  By the inductive hypothesis there exists a basis $v_1,\ldots,v_{d-1}$ of $W_2$ containing a basis of each of $W_3,\ldots, W_h$ and $W_1',\ldots,W_h'$.  If $W_i^*\subset W_2$ for $i>1$ then $W_i'=W_i^*$ for $i>1$.  So in this case if we complete $v_1,\ldots,v_{d-1}$ to any basis of $V$ we are done.  Otherwise, let $l$ be the maximal index with $W_l^* \not\subset W_2$ and let $v_d\in W_l^*\backslash W_l'$.  We claim that $B=\{v_1,\ldots,v_d\}$ is a basis of $V$ with the required property.  It clearly contains a basis of $W_i$ for each $i$.  Let $i\in \{1,\ldots,h^*\}$.  If $i>l$ then $W_i^*=W_i'$ and so by construction $B$ contains a basis of $W_i^*$.  If $i\leq l$ then $v_d\in W_l^*\backslash W_l' \subset W_i^*\backslash W_i'$.  Since $B$ contains a basis $B_i'$ of $W_i'$ and $W_i'$ is a hyperplane in $W_i^*$, we see that $B_i'\cup \{v_d\}$ is a basis of $W_i^*$. 
\end{proof}
Using our notation from the last section, suppose that for $P\in D$ we have $D_P=D_{P,1}+D_{P,2}$ where $D_{P,1}$ and $D_{P,2}$ are effective divisors with no irreducible components in common.  We may then prove the following versions of Lemma \ref{large} and Theorem \ref{cor2}.
\begin{lemma}
\label{large2}
Let $D=\sum_{i=1}^r D_i$ be a nonzero divisor on a nonsingular variety $X$ with $D_i$ effective for each $i$.  Let $P\in D$.  Let $f_{P,j}(m,n)=l(nD-mD_{P,j})-l(nD-(m+1)D_{P,j})$ for $j=1,2$.  If there exists $n>0$ such that either $\sum_{m=0}^{\infty}(m-n)f_{P,j}(m,n)>0$ or $D_{P,j}=0$ for all $P\in D$ and $j=1,2$ then $nD$ is very large.
\end{lemma}
\begin{theorem}
\label{cor22}
Let $X$ be a nonsingular variety.  Let $q= \dim X$.  Let $D=\sum_{i=1}^{r}D_i$ be a divisor on $X$ such that $D_{P,j}$ is nef for all $P\in D$ and $j=1,2$.  Suppose also that every irreducible component of $D$ is nonsingular.  If
\begin{equation*}
D^q>2q D^{q-1}.D_{P,j}, \qquad \forall P\in D, j=1,2
\end{equation*}
then $nD$ is very large for $n\gg 0$.
\end{theorem}
The proofs are similar to the proofs of Lemma \ref{large} and Theorem \ref{cor2}.  The only difference is that in the proof of Lemma \ref{large2}, we look at the two filtrations of $L(nD)$ given by $W_j=L(nD-jD_{P,1})$ and $W_j^*=L(nD-jD_{P,2})$ and we use the filtration lemma to construct a basis $f_1,\ldots,f_{l(nD)}$ that contains a basis for each $W_j$ and $W_j^*$.

Suppose now that $D=\sum_{i=1}^rD_i$ where the $D_i$'s are effective divisors and the intersection of any $m+1$ distinct $D_i$'s is empty.  We may then write $D_P=D_{P,1}+D_{P,2}$ where $D_{P,1}$ and $D_{P,2}$ are each not a sum of more than $[\frac{m+1}{2}]$ $D_i$'s, where $[x]$ denotes the greatest integer in $x$.  Using this, we get the following improvements to the part (a)'s of Theorems \ref{cor4} and \ref{cor4b}.
\renewcommand{\thetheorem}{\arabic{section}.\arabic{theorem}A}
\setcounter{theoremb}{\value{theorem}}
\begin{theorema}
\label{cor42}
Let $X$ be a projective variety defined over a number field $k$.  Let $q=\dim X$.  Let $D=\sum_{i=1}^{r} D_i$ be a divisor on $X$ defined over $k$ such that the $D_i$'s are effective Cartier divisors with no irreducible components in common and the intersection of any $m+1$ distinct $D_i$'s is empty.  If $D_i$ is quasi-ample for each $i$ and $r> 2[\frac{m+1}{2}]q$ then $X\backslash D$ is quasi-Mordellic.
\end{theorema}
\begin{theoremb}
\label{cor4b2}
Let $X$ be a complex projective variety.  Let $q=\dim X$.  Let $D=\sum_{i=1}^{r} D_i$ be a divisor on $X$ such that the $D_i$'s are effective Cartier divisors with no irreducible components in common and the intersection of any $m+1$ distinct $D_i$'s is empty.  If $D_i$ is quasi-ample for each $i$ and $r>2[\frac{m+1}{2}]q$ then $X \backslash D$ is quasi-Brody hyperbolic.
\end{theoremb}
Unfortunately, we need the requirement that the $D_i$'s have no irreducible components in common so that we may have $D_{P,1}$ and $D_{P,2}$ with no irreducible components in common (which is necessary in proving Lemma \ref{large2}).  Because of this, we cannot prove a finiteness result about ample divisors as we did in the last section, since the restrictions of the $D_i$'s to a subvariety of $X$ may have irreducible components in common.

\section{Surfaces}
\label{ssurf}
\renewcommand{\thetheorem}{\arabic{section}.\arabic{theorem}}
We will now see that we may make the results of the last two sections more precise if we restrict to the case where $X$ is a surface.  With regards to integral points, this section builds on some of the work in \cite{Co2}.  Corvaja and Zannier prove, essentially, Theorem \ref{surf} \cite[Main Theorem]{Co2} and they prove Theorem \ref{surf3a} when $m=2$ and the $D_i$'s have multiples which are all numerically equivalent.  The Nevanlinna-theoretic analogues of the results in \cite{Co2} were proved by Ru and Liu in \cite{Ru4}.  Our results overlap with their results as well.

We first prove a consequence of the Hodge Index theorem.
\begin{lemma}
\label{Hodge}
Let $D$ be a divisor on a nonsingular surface $X$ with $D^2>0$.  Then $(D^2)(E^2)\leq (D.E)^2$ for any divisor $E$ on $X$.
\end{lemma}
\begin{proof}
By the Hodge index theorem, the intersection pairing on Num $X\bigotimes\mathbb{R}$ can be diagonalized with one $+1$ on the diagonal and all other diagonal entries $-1$.  We will identify elements of Pic $X$ as elements of Num $X\bigotimes\mathbb{R}$ in the canonical way.  Extend $D$ to an orthogonal basis $B$ of Num $X\bigotimes\mathbb{R}$.  Let $E$ be any divisor on $X$.  Writing $E$ in the basis $B$, it is apparent from the Hodge index theorem that $(D^2)(E^2)\leq (D.E)^2$.
\end{proof}

For surfaces, the more precise version of Theorem~\ref{cor22} is
\begin{theorem}[Corvaja-Zannier]
\label{surf}
Let $X$ be a nonsingular projective surface.  Let $D=\sum_{i=1}^{r}D_i$ be a nef divisor on $X$ with the $D_i$'s effective divisors and $D^2>0$.  For $P\in D$, let $D_P=\sum_{i:P\in D_i}D_i=D_{P,1}+D_{P,2}$ where $D_{P,1}$ and $D_{P,2}$ are effective divisors with no irreducible components in common.  Suppose that for all $P\in D$, $j=1,2$ and $m,n>0$ we have either $l(nD-mD_{P,j})=0$ or
\begin{equation*}
l(nD-mD_{P,j})-l(nD-(m+1)D_{P,j})\leq (nD-mD_{P,j}).mD_{P,j}+O(1)
\end{equation*}
where the constant does not depend on $m$ or $n$.  Let $A_{P,j}=D_{P,j}.D_{P,j},B_{P,j}=D.D_{P,j}$, and $C=D.D$ for $j=1,2$.  If for all $P \in D$ and $j=1,2$ either we have $D_{P,j}=0$ or we have 
\begin{align*}
&A_{P,j}>0 \Longrightarrow B_{P,j}^2-2A_{P,j}C+3A_{P,j}B_{P,j}+(3A_{P,j}-B_{P,j})\sqrt{B_{P,j}^2-A_{P,j}C}<0\\
&A_{P,j}=0 \Longrightarrow C>4B_{P,j}\\
&A_{P,j}<0 \Longrightarrow B_{P,j}^2-2A_{P,j}C+3A_{P,j}B_{P,j}+(3A_{P,j}-B_{P,j})\sqrt{B_{P,j}^2-A_{P,j}C}>0
\end{align*}
then $nD$ is very large for $n\gg 0$ (note that by Lemma \ref{Hodge} $B_{P,j}^2-A_{P,j}C>0$).
\end{theorem}
\begin{proof}
Let $P \in D$ and $j\in \{1,2\}$ with $D_{P,j}\neq 0$.  Let $A=A_{P,j}$ and $B=B_{P,j}$.  By assumption, we have
\begin{align*}
f_{P,j}(m,n)&=\dim H^0(X,\co(nD-m D_{P,j}))-\dim H^0(X,\co(nD-(m+1)D_{P,j}))\\
&\leq nB-mA+O(1)
\end{align*}
where the constant in the $O(1)$ does not depend on $m$ or $n$.  We have
\begin{equation*}
l(nD)=\frac{D^2}{2}n^2+O(n)=\frac{C}{2}n^2+O(n).
\end{equation*}
Solving
\begin{equation*}
\sum_{m=0}^{M(n)} nB-mA+O(1)=\frac{C}{2}n^2+O(n)= l(nD)
\end{equation*}
for $M(n)$, we get
\begin{align*}
&M(n)=\frac{B \pm \sqrt{B^2-AC}}{A}n+O(1), &A \neq 0\\
&M(n)=\frac{C}{2B}n+O(1),  &A=0,B\neq 0\\
&M(n)=O(n^2),  &A=0,B=0.
\end{align*}
From now on, we will always choose the minus sign in the first expression above.  We also have $\sum_{m=0}^{\infty}f_{P,j}(m,n)=l(nD)$.  Therefore by Lemma \ref{CZlemma},
\begin{equation}
\label{MP}
\sum_{m=0}^{\infty}(m-n)f_{P,j}(m,n) \geq \sum_{m=0}^{M(n)}m(nB-mA+O(1))-nl(nD).
\end{equation}
Let $K=\frac{B - \sqrt{B^2-AC}}{A}$.  If $A \neq 0$ then substituting $K$ into (\ref{MP}) we get
\begin{equation*}
\sum_{m=0}^{\infty}(m-n)f_{P,j}(m,n)\geq(-\frac{A}{3}K^3+\frac{B}{2}K^2-\frac{C}{2})n^3+O(n^2)
\end{equation*}
So if $-\frac{A}{3}K^3+\frac{B}{2}K^2-\frac{C}{2}>0$ then by Lemma \ref{large2}, $nD$ will be very large for $n\gg 0$.  Algebraic simplification then gives the theorem in the case $A \neq 0$.  The other cases are similar.
\end{proof}
\begin{lemma}
\label{clemma}
Let $X$ be a nonsingular projective surface.  Let $C$ be an irreducible curve on $X$ and $D$ any divisor on $X$.  Then
\begin{equation*}
h^0(D)-h^0(D-C)\leq \max\{0,1+C.D\}.
\end{equation*}
\end{lemma}
\begin{proof}
The statement depends only on the linear equivalence class of $D$, so replacing $D$ by an appropriate divisor linearly equivalent to $D$, we may assume that the support of $D$ does not contain any possible singularity of $C$.  By Lemma \ref{exact} we have
\begin{equation*}
h^0(D)-h^0(D-C)\leq \dim H^0(C,\co(D)|_C)
\end{equation*}
Since the support of $D$ does not contain any singularity of $C$, $\co(D)|_C$ has degree $C.D$ on $C$ and $\dim H^0(C,\co(D)|_C)\leq \max\{0,1+C.D\}$.
\end{proof}
\begin{lemma}
\label{slemma}
Let $X$ be a nonsingular projective surface.  Let $D$ be a nef divisor on $X$.  Let $E$ be an effective divisor on $X$ such that either $E$ is linearly equivalent to an irreducible curve or for every irreducible component $C$ of $E$, $C.E\leq 0$.  Then for all $m,n>0$ either $l(nD-mE)=0$ or
\begin{equation}
\label{ineq}
l(nD-mE)-l(nD-(m+1)E)\leq (nD-mE).E+O(1)
\end{equation}
where the constant is independent of $m$ and $n$.
\end{lemma}
\begin{proof}
In the first case, suppose $E$ is linearly equivalent to an irreducible curve $C$.  If $(nD-mE).E\geq 0$ then (\ref{ineq}) holds by Lemma \ref{clemma}.  If $(nD-mE).E=nD.C-mC.C<0$ then since $D$ is nef, we must have $C.C>0$.  But if $l(nD-mE)>0$ then $nD-mE$ is linearly equivalent to an effective divisor $F=G+mC$ where $m\geq 0$ and $G$ is an effective divisor not containing $C$.  Since clearly $G.C\geq 0$, $F.C=(nD-mE).E<0$ implies $C.C<0$, a contradiction.  So either $l(nD-mE)=0$ or (\ref{ineq}) holds in this case.

Now suppose we are in the second case, where for every irreducible component $C$ of $E$, $C.E\leq 0$.  Let $E=\sum_{j=1}^{k}a_j C_j$, where each $C_j$ is a distinct prime divisor.  Then as in the proof of Theorem~\ref{cor2} we have
\begin{multline*}
l(nD-mE)-l(nD-(m+1)E) \leq \\
\sum_{j=1}^k \sum_{l=0}^{a_{j}-1}\dim H^0(C_{j},i^*_{C_{j}}\co(nD-mE -\sum_{j'=1}^{j-1}a_{j'}C_{j'}-lC_{j}))
\end{multline*}
But
\begin{align*}
\dim H^0(C_{j},i^*_{C_{j}}\co(nD-mE -\sum_{j'=1}^{j-1}a_{j'}C_{j'}-lC_{j}))&\leq \dim H^0(C_{j},i^*_{C_{j}}\co(nD-mE))+O(1)\\
&\leq (nD-mE).C_j+O(1),
\end{align*}
where the constant is independent of $m$ and $n$.  The second inequality follows since $(nD-mE).C_j\geq nD.C_j\geq 0$ as $D$ is nef and $E.C_j\leq 0$.  Combining the above inequalities, we then see that (\ref{ineq}) always holds in this case.
\end{proof}
Going back to the General Setup of Section \ref{gsetup}, we have
\renewcommand{\thetheorem}{\arabic{section}.\arabic{theorem}A}
\setcounter{theoremb}{\value{theorem}}
\begin{theorema}
\label{surf3a}
Let $X$ be a projective surface.  Suppose the $D_i$'s have no irreducible components in common.\\\\
(a). If $D_i$ is quasi-ample for all $i$ and $r\geq 4[\frac{m+1}{2}]$ then $X\backslash D$ is quasi-Mordellic.\\
(b). If $D_i$ is ample for all $i$ and either $m$ is even and $r>2m$ or $m$ is odd and $r>2m+1$ then $X\backslash D$ is Mordellic.
\end{theorema}
\begin{theoremb}
\label{surf3b}
Let $X$ be a projective surface.  Suppose the $D_i$'s have no irreducible components in common.\\\\
(a). If $D_i$ is quasi-ample for all $i$ and $r\geq 4[\frac{m+1}{2}]$ then $X\backslash D$ is quasi-Brody hyperbolic.\\
(b). If $D_i$ is ample for all $i$ and either $m$ is even and $r>2m$ or $m$ is odd and $r>2m+1$ then $X\backslash D$ is complete hyperbolic and hyperbolically imbedded in $X$.  In particular, $X\backslash D$ is Brody hyperbolic.
\end{theoremb}
\begin{proof}
We'll prove the part (a)'s first.  It suffices to prove these in the case $r=4[\frac{m+1}{2}]$.  As in the proofs of Theorems \ref{cor4} and \ref{cor4b}, we may use Lemma \ref{reduce} to reduce to the case where $X$ is nonsingular, $|D_i|$ is base-point free for all $i$, and $\dim \Phi_{D_i}(X)=2$ for all $i$.  Therefore $D_i^2>0$ and $D_i$ is nef for each $i$.  By Lemma \ref{equi}, $D$ is equidegreelizable.  So we may find positive integers $a_1,\ldots,a_r$ such that if $D'=\sum_{i=1}^{r}a_iD_i$ then $\frac{a_iD_i.D'}{(D')^2}$ is arbitrarily close to $\frac{1}{r}$ for all $i$.  Since at most $m$ $D_i$'s meet at any given point, $D'_P$ is a sum of at most $m$ $a_iD_i$'s for any $P\in D'$.  Therefore we may write $D'_P=D'_{P,1}+D'_{P,2}$ where each $D'_{P,j}$ is a sum of at most $[\frac{m+1}{2}]$ $a_iD_i$'s, and $D'_{P,1}$ and $D'_{P,2}$ have no irreducible components in common.  Note that when $D'_{P,j}\neq 0$, we have, from our assumptions on the $D_i$'s, that $|D'_{P,j}|$ is base-point free and $\dim \Phi_{D'_{P,j}}(X)=2$.  So by Theorem \ref{Bertini}, $D'_{P,j}$ is linearly equivalent to an irreducible curve.  Therefore, by Lemma \ref{slemma}, we will be able to apply Theorem \ref{surf} to $D'$.

The hardest case is clearly when $D'_{P,j}$ is a sum of the maximum $[\frac{m+1}{2}]$ $a_iD_i$'s.  For simplicity, we will now restrict to this case.  It follows that in the notation of Theorem \ref{surf} we may take, for all such $P$ and $j$,
\begin{equation*}
\left|\frac{C}{B_{P,j}}-\frac{r}{[\frac{m+1}{2}]}\right|=\left|\frac{C}{B_{P,j}}-4\right|<\epsilon
\end{equation*}
where by adjusting the $a_i$'s in $D'$, $\epsilon$ may be made arbitrary close to $0$ while at the same time $\frac{A_{P,j}}{B_{P,j}}$ is positive and bounded away from $0$.  Furthermore, by Lemma \ref{Hodge}, $\frac{A_{P,j}}{B_{P,j}}\leq \frac{B_{P,j}}{C}$.  Let $a=\frac{A_{P,j}}{B_{P,j}}$ and $c=\frac{C}{B_{P,j}}$.  Then by Theorem \ref{surf}, we must show that
\begin{equation}
\label{sineq}
1-2ac+3a+(3a-1)\sqrt{1-ac}<0
\end{equation}
where $0<a\leq\frac{1}{c}$.  When $c=4$, we get $1-5a+(3a-1)\sqrt{1-4a}$, which is easily seen to have a root only at $a=0$ for $0\leq a\leq\frac{1}{4}$, and is negative for $0<a\leq\frac{1}{4}$ since putting $a=\frac{1}{4}$ gives $-\frac{1}{4}$.  So when $c=4+\epsilon$, since $a$ is bounded away from zero as $\epsilon \to 0$, we see that (\ref{sineq}) is negative for small enough $\epsilon$.  Therefore by Theorem \ref{surf}, $nD'$ is very large for $n\gg 0$.  Since $D'$ is quasi-ample, $\Phi_{nD'}$ is a birational map to projective space for some arbitrarily large $n$.  Therefore by Theorems \ref{maina} and \ref{mainb} we are done, as $D$ and $D'$ have the same support.

Assume the hypotheses in the part (b)'s.  Let $Y$ be the Zariski-closure of a set of $D$-integral points (resp. $f(\mathbb{C})$).  By what we have proven above, $\dim Y\leq 1$.  If $\dim Y=1$, let $C$ be an irreducible component of this curve with $\dim C>0$.  Since each $D_i$ is ample, $D_i$ must intersect $C$ in a point.  Since at most $m$ $D_i$'s meet at a point and $r>2m$, we see that $D|_C$ contains at least $3$ distinct points.  Therefore by Siegel's (resp. Picard's) theorem we get a contradiction as the above gives a dense set of $D|_C$-integral points (resp. a dense holomorphic map $\mathbb{C} \to C\backslash D|_C$).  This same argument and Theorem \ref{hyperbolic} show that in the analytic case $X\backslash D$ is hyperbolic and hyperbolically embedded in $X$.
\end{proof}
It is possible to make minor improvements to this theorem.  For example,
\begin{theorema}
\label{surf4a}
Let $X$ be a nonsingular projective surface.  Suppose $m=2$, $D=\sum_{i=1}^4D_i$, $D_i.D_j>0$ for $i\neq j$, $D_1^2>0$, $D_i$ is nef for all $i$, and the $D_i$'s have no irreducible components in common.  Suppose also that the conclusion of Lemma \ref{slemma} holds with $D$ any positive integral linear combination of the $D_i$'s and $E=D_i$, for $i=1,2,3,4$.  Then $X\backslash D$ is quasi-Mordellic.
\end{theorema}
\begin{theoremb}
\label{surf4b}
With the same hypotheses as above, in the analytic setting, $X\backslash D$ is quasi-Brody hyperbolic.  
\end{theoremb}
\begin{proof}
We first show that for any $\epsilon>0$, $(\sum_{i=1}^4e^{a_i}D_i)^2\geq e^{\frac{2}{3}\max_i\{a_i\}}$ on the plane $(1+\epsilon)a_1+\sum_{i=2}^4a_i=0$.  If $\max_i\{a_i\}=a_1$ then $(\sum_{i=1}^4e^{a_i}D_i)^2\geq e^{2a_1}D_1^2\geq e^{2a_1}$.  Otherwise, if $\max_i\{a_i\}=a_j$, $j>1$, then clearly we must have $a_k\geq -\frac{a_j}{3}$ for some $j\neq k$.  Then $(\sum_{i=1}^4e^{a_i}D_i)^2\geq e^{a_j+a_k}D_j.D_k\geq e^{\frac{2}{3}a_j}$ since $D_j.D_k\geq 1$.  Therefore $(\sum_{i=1}^4e^{a_i}D_i)^2$ takes a minimum on the plane $(1+\epsilon)a_1+\sum_{i=2}^4a_i=0$.  So looking at the Lagrange multiplier equations as in Lemma \ref{equi}, there exist real numbers $b_i>0,\lambda>0$ (depending on $\epsilon$) such that if $D'=\sum_{i=1}^4b_iD_i$ then $b_1D_1.D'=(1+\epsilon)\lambda$ and $b_iD_i.D'=\lambda$ for $i=2,3,4$, or written differently, $\frac{(D')^2}{b_1D_1.D'}=\frac{4+\epsilon}{1+\epsilon}$ and $\frac{(D')^2}{b_iD_i.D'}=4+\epsilon>4$ for $i=2,3,4$.  Note also that it follows from the inequality we proved above that in terms of $a_1,\ldots,a_4$, the region where $(\sum_{i=1}^4e^{a_i}D_i)^2$ takes a minimum may be bounded independently of $\epsilon$.  Therefore there exist positive constants $K,K'$ independent of $\epsilon$, such that we may choose $K<b_i<K'$ for all $i$, and in particular, as $\epsilon\to 0$, $\frac{(b_1D_1)^2}{b_1D_1.D'}$ is bounded away from zero.

We now choose positive integers $c_i$ such that $\frac{c_i}{c_j}$ closely approximates $\frac{b_i}{b_j}$, and let $E=\sum_{i=1}^4c_iD_i$.  Having chosen $\epsilon$ small enough and the integers $c_i$ correctly, we will then have $E^2>4c_iD_i.E$ for $i=2,3,4$ and we will have $\frac{E^2}{c_1D_1.E}$ close enough to $4$ (see the proof of Theorem \ref{surf3a}, B) so that the inequalities in $\ref{surf}$ hold for $E_{P,j}=c_iD_i$ for any $i$.  Since $m=2$, we may always take $E_{P,j}=0$ or $E_{P,j}=c_iD_i$ for some $i$.  By our hypotheses, we may apply Theorem \ref{surf}, so $nE$ is very large for $n\gg 0$.  Since $D_1^2>0$, $E$ is quasi-ample.  So we are done by Theorems \ref{maina} and \ref{mainb}, as $D$ and $E$ have the same support.
\end{proof}
\renewcommand{\thetheorem}{\arabic{section}.\arabic{theorem}}
\begin{example}
Let $X=\mathbb{P}^1\times \mathbb{P}^1$.  Let $D_1=\{0\}\times \mathbb{P}^1, D_2=\mathbb{P}^1\times \{0\}$, and let $D_3$ and $D_4$ be ample effective divisors on $X$.  Suppose also that the intersection of any three of the $D_i$'s is empty.  Let $D=\sum_{i=1}^4D_i$.  Then the hypotheses of Theorems \ref{surf4a}, B are satisfied and $X\backslash D$ is quasi-Mordellic and quasi-Brody hyperbolic.  Note also that $X\backslash D_1\cup D_2\cong \mathbb{A}^2\cong \mathbb{P}^2\backslash \{\text{a line}\}$.  Therefore, we can also prove many theorems for $\mathbb{P}^2\backslash D$, where $D$ is a sum of three effective divisors on $\mathbb{P}^2$.
\end{example}
Recently, Corvaja and Zannier \cite{Co6} have shown another way their methods may get results on $\mathbb{P}^2\backslash D$ where $D$ is a sum of three effective divisors satisfying certain hypotheses.

We have the following general corollary to the above theorems.
\renewcommand{\thetheorem}{\arabic{section}.\arabic{theorem}A}
\setcounter{theoremb}{\value{theorem}}
\begin{corollarya}
Let $X$ be a projective surface.  Suppose $m=2$, $D=\sum_{i=1}^4D_i$, $D_1,D_2,D_3$ are quasi-ample, $\kappa(D_4)>0$, and the $D_i$'s have no irreducible components in common.  Then $X\backslash D$ is quasi-Mordellic.
\end{corollarya}
\begin{corollaryb}
Let $X$ be a projective surface.  Suppose $m=2$, $D=\sum_{i=1}^4D_i$, $D_1,D_2,D_3$ are quasi-ample, $\kappa(D_4)>0$, and the $D_i$'s have no irreducible components in common.  Then $X\backslash D$ is quasi-Brody hyperbolic.
\end{corollaryb}
\begin{proof}
We first reduce to the situation of Lemma \ref{reduce}.  So $X$ is nonsingular, each $D_i$ is nef, and $D_1^2,D_2^2,D_3^2>0, D_4^2\geq 0$  By Lemma \ref{m1}, $D_i.D_j>0$ for $i\neq j$.  For $i=1,2,3$ and $n>0$ $nD_i$ is linearly equivalent to an irreducible curve by Theorem \ref{Bertini}, since by our reductions $|nD_i|$ is base-point free and $\dim \Phi_{nD_i}(X)=2$.  The same holds for $nD_4$ if $D_4^2>0$.  If $D_4^2=0$, then for every irreducible component $C$ of $D_4$ we must have $C.D_4=0$ since $D_4$ is nef.  This verifies the hypotheses of Lemma \ref{slemma} with $E=D_i$ for $i=1,2,3,4$.  Therefore, we may apply Theorems \ref{surf4a}, B.
\end{proof}
We note that one can construct examples where $m=2$, $D_1$ and $D_2$ are quasi-ample, $\kappa(D_3)=\kappa(D_4)=1$, the $D_i$'s have no irreducible components in common, and there exist dense sets of $D$-integral points.  We now prove a theorem in the case where we only have $\kappa(D_i)>0$ for all $i$.
\begin{theorema}
Let $X$ be a projective surface.  Suppose the $D_i$'s have no irreducible components in common.  If $\kappa(D_i)>0$ for all $i$ and $r> 4[\frac{m+1}{2}]$ then there does not exist a Zariski-dense set of $D$-integral points on $X$.
\end{theorema}
\begin{theoremb}
Let $X$ be a projective surface.  Suppose the $D_i$'s have no irreducible components in common.  If $\kappa(D_i)>0$ for all $i$ and $r> 4[\frac{m+1}{2}]$ then there does not exist a holomorphic map $f:\mathbb{C}\to X\backslash D$ with Zariski-dense image.
\end{theoremb}
\begin{proof}
We first reduce to the situation of Lemma \ref{reduce}, so in particular $|D_i|$ is base-point free for all $i$.  In this case, for any subset $I\subset \{1,\ldots,r\}$, if $D_I=\sum_{i\in I}D_i$ is quasi-ample, there exists $n_I>0$ such that $\dim \Phi_{n_ID_I}(X)=2$.  Since the $D_i$'s are nef, this happens if and only if $D_i.D_j>0$ for some $i,j\in I$.  Let $N=\prod_{I}n_I$ where the $I$ ranges over subsets such that $D_I$ is quasi-ample.  Let $D'=ND$ and $D_i'=ND_i$.  Then we see that for any nonnegative integral linear combination, $E$, of the $D_i'$'s, if $E$ is quasi-ample, then $E$ is linearly equivalent to an irreducible divisor since $|E|$ is base-point free and $\dim \Phi_E(X)=2$, and otherwise, for every irreducible component $C$ of $E$ we have $C.E=0$.  Therefore, by Lemma \ref{slemma}, replacing $D$ by $D'$, we may assume that we may apply Theorem \ref{surf} to any nonnegative linear combination of the $D_i$'s.

By Theorem \ref{tm1}, we are done if any three of the $D_i$'s have pairwise empty intersection.  So suppose that this is not the case.  Then we have $m\geq 2$ and $r\geq 5$.  We now show that $D$ is equidegreelizable.  As in the proof of Lemma \ref{equi}, it suffices to show that $(\sum_{i=1}^re^{a_i}D_i)^2$ attains a minimum on the plane $\sum_{i=1}^ra_i=0$.  For this, it will suffice to show that $(\sum_{i=1}^re^{a_i}D_i)^2\geq e^{\frac{1}{3}\max_i\{a_i\}}$.  Suppose $\max_i\{a_i\}=a_j$ for $j\in \{1,\ldots,r\}$.  Let $a_k$ and $a_l$ be some choice of the next largest $a_i$'s.  Clearly, since $\sum_{i=1}^ra_i=0$, we must have $a_k,a_l\geq -\frac{2a_j}{r-2}\geq -\frac{2}{3}a_j$ since $r\geq 5$.  We now show that either $D_j.D_k\geq 1$ or $D_j.D_l\geq 1$.  Suppose otherwise.  Then by our assumption, we must have $D_k.D_l\geq 1$.  But then $D_k+D_l$ is quasi-ample, and so we must have $(D_k+D_l).D_j\geq 1$ by Lemma \ref{m1}, a contradiction.  So if, say, $D_j.D_k\geq 1$ then $(\sum_{i=1}^re^{a_i}D_i)^2\geq e^{a_j+a_k}D_j.D_k\geq e^{\frac{1}{3}\max_i\{a_i\}}$, as was to be shown.  Since $D$ is equidegreelizable, there exist positive integers $c_i$ such that $D'=\sum_{i=1}^rc_iD_i$, and $\frac{c_iD_i.D'}{(D')^2}$ is as close as we like to $\frac{1}{r}$.  Since we may choose $D'_{P,j}$ to consist of a sum of at most $[\frac{m+1}{2}]$ $c_iD_i$'s and $r>4[\frac{m+1}{2}]$, we may choose the $c_i$'s so that we always have $C>4B_{P,j}$.  We also have $A_{P,j}\geq 0$.  But then, as we have seen previously, the inequalities of Theorem \ref{surf} will be satisfied.
\end{proof}
\renewcommand{\thetheorem}{\arabic{section}.\arabic{theorem}}
To summarize some of the results in this section:
\begin{theorem}
Let $X$ be a projective surface.  Suppose that $m\leq 2$ and the $D_i$'s have no irreducible components in common.  Then all of the Main Conjectures (Conjectures \ref{conjmaina},B-\ref{conj2a},B) are true.
\end{theorem}
\section{Small S}
\label{ssmall}
We now prove some theorems in the special case that $\#S$ is small relative to the number of components of $D$.  Throughout we use the general Diophantine setup of Section \ref{gsetup}.
\begin{theorem}
\label{sthm}
Suppose that $D_i$ is defined over $k$ for all $i$.  Let $s=\#S$.\\\\
(a).  If $D_i$ is quasi-ample for all $i$ and $r>ms$ then there exists a proper closed subset $Z\subset X$ such that for any set $R$ of $(D,S)$-integral points on $X$, $R\backslash Z$ if finite.\\
(b).  If $D_i$ is ample for all $i$ and $r>ms$ then all sets of $(D,S)$-integral points on $X$ are finite.
\end{theorem}
\begin{proof}
We reduce to the case where $X$ is nonsingular.  We prove part (a) first.  Our proof is a modification of the proof of Theorem~\ref{maina}.  Suppose $R$ is a Zariski-dense set of $(D,S)$-integral points on $X$.  Then as in the proof of Theorem \ref{maina}, there exists a sequence $P_i$ in $R$ such that for each $v$ in $S$, $\{P_i\}$ converges to a point $P_v\in X(k_v)$ and $\bigcup\{P_i\}$ is Zariski-dense in $X$.  Since $r>ms$, there exists an index $i$ such that $P_v\notin D_i(k_v)$ for all $v\in S$.  Since $D_i$ is quasi-ample, it follows from Lemma \ref{exact} and the argument in Lemma \ref{nefbig} that for some $n>0$, $l(nD_i-\sum_{j\neq i}D_i)>0$.  Then $(n+1)D_i-\sum_{j\neq i}D_i$ is quasi-ample, and so for some $n'>0$, and $E=n'(n+1)D_i-n'\sum_{j\neq i}D_i$, $\Phi_{E}$ is birational.  Now let $S'$ be the set of places $v\in S$ such that $P_v \in D(k_v)$.  Let $\phi_1,\ldots,\phi_{l(E)}$ be a basis for $L(E)$ over $k$.  Then for any $v\in S'$, $\text{ord}_F\prod_{i=1}^{l(E)}\phi_i>0$ for every irreducible component $F$ of $D$ such that $P_v \in F(k_v)$.  This is precisely what we used the largeness hypothesis for in the proof of Theorem \ref{maina}.  Let $\phi=(\phi_1,\ldots,\phi_{l(E)})$.  Let $L_{jv}=\phi_j$ for $j=1,\ldots,l(E)$ and $v\in S$.  Then the same proof as in Theorem \ref{maina} (replacing $D$ by $E$ in appropriate places) proves part (a).

For (b), let $R$ be a set of $(D,S)$-integral points on $X$.  Let $Y$ be an irreducible component of the Zariski-closure, $\overline{R}$, of $R$.  Suppose $\dim Y>0$.  Then $D$ pulls back to a sum of $r$ ample effective divisors on $Y$ such that at most $m$ of them meet at a point.  But then part (a) applied to $D|_Y$ contraticts the fact that $R\cap Y$ is a dense set of $(D|_Y,S)$-integral points.  Therefore $\dim Y=0$.
\end{proof}
When $\#S=1$ this theorem gives a particularly strong result.
\begin{corollary}
\label{qs1}
Suppose $\#S=1$.  If $D_i$ is ample for all $i$ and $r>m$ then all sets of $(D,S)$-integral points on $X$ are finite.
\end{corollary}
It follows from the Dirichlet unit theorem that $\#S=1$ if and only if $\co_{k,S}^*$ is finite if and only if $\co_{k,S}=\mathbb{Z}$ or the ring of integers of a complex quadratic field.  The inequality in Corollary \ref{qs1} is sharp as the next example shows.
\begin{example}
Let $X=\mathbb{P}^n$.  Let $k=\mathbb{Q}$ and let $S$ consist only of the prime at infinity.  Let $D_i$ be the divisor on $\mathbb{P}^n$ defined by $x_i=0$, where $x_0,\ldots, x_n$ are homogeneous coordinates on $\mathbb{P}^n$.  Let $D=\sum_{i=1}^n a_iD_i$.  Let $m=\sum_{i=1}^na_i$.  Then the set of points with $x_0\in \mathbb{Z}$ and $x_i=1$, $i=1,\ldots,n$, is an infinite set of $(D,S)$-integral points on $X$ and $D$ is a sum of $m$ ample divisors defined over $\mathbb{Q}$.
\end{example}
\begin{theorem}
\label{Pic}
Let $X$ be a nonsingular projective variety.  Suppose $\#S=1$.  Let $\rho$ denote the Picard number of $X$ and let $n$ be the rank of the group of $k$-rational points of $\text{Pic}^0(X)$.  Suppose that the $D_i$'s are defined over $k$ for all $i$ and have no irreducible components in common.  If $r>\rho+n$ then there does not exist a dense set of $(D,S)$-integral points on $X$.
\end{theorem}
Our proof is essentially the first half of the proof of Theorem 2.4.1 in \cite{Vo2}.
\begin{proof}
It follows from the definitions that the group of divisor classes with a representative defined over $k$ has rank at most $\rho+n$.  Since $r>\rho+n$, there exists a linear combination of the $D_i$'s that is principal, equal to $(f)$ for some nonconstant rational function $f$ on $X$.  Let $R$ be a set of $(D,S)$-integral points on $X$.  Since all of the poles of $f$ lie in $D$ there exists an $a\in k$ such that $af$ takes on integral values on $R$.  Since the poles of $\frac{1}{f}$ also lie in $D$, the same reasoning applies to $\frac{1}{f}$.  Therefore $f(R)$ lies in only finitely many cosets of the group of units $\co_{k,S}^*$.  But since $\#S=1$, $\co_{k,S}^*$ is finite.  Therefore $R$ lies in the finite union of proper subvarieties of $X$ of the form $f=a$ for a finite number of $a\in k$.
\end{proof}
We note that the requirement in all of these results that not only $D$ be defined over $k$, but that the $D_i$'s be defined over $k$ is absolutely necessary.  For example, if $X=\mathbb{P}^1$, $k=\mathbb{Q}$, $S=\{\infty\}$, and $D=P+Q$ where $P$ and $Q$ are conjugate over a real quadratic field, then from Pell's equation there do exist dense sets of $(D,S)$-integral points on $X$.
\section{Results on the General Conjectures}
\label{SVgeneral}
We will now consider the case where the integral points are allowed to vary over number fields of a bounded degree over some number field $k$.  As an application of their results on surfaces in \cite{Co2}, Corvaja and Zannier prove
\begin{theorem}
Let $X$ be a projective curve defined over a number field $k$.  Let $S$ be a finite set of places of $k$ containing the archimedean places.  Let $D=\sum_{i=1}^r P_i$ be a divisor on $X$ defined over $k$ such that the $P_i$'s are distinct points.  If $r>4$ then all sets of $D$-integral points on $X$ quadratic over $k$ are finite.
\end{theorem}
This theorem can also be obtained as a consequence of a result of Vojta (see Section \ref{sVo}).  Using the same technique Corvaja and Zannier used, looking at symmetric powers of $X$, our higher-dimensional results give
\begin{theorem}
Let $n=\dim X$.  If $D_i$ is ample for all $i$ and $r>2d^2mn$ then all sets of $D$-integral points on $X$ of degree $d$ over $k$ are finite.
\end{theorem}
\begin{proof}
Suppose $r>2d^2mn$ and let $R\subset X(\kk)$ be a set of $D$-integral points on $X$ of degree $d$ over $k$.  It suffices to prove the finiteness of $R$ in the case where for every $P\in R$ we have $[k(P):k]=d$.  Let $X^d$ be the $d$-fold product of $X$ with itself, and let $\pi_i:X^d\to X$ be the $i$-th projection map for $i=1,\ldots, d$.  Let $\Sym^d X$ denote the $d$-fold symmetric product of $X$ with itself and let $\phi:X^d\to \Sym^d X$ be the natural map.  Let $E_i=\phi(\pi_1^*D_i)$ and $E=\sum_{i=1}^rE_i$.  We have that $\phi^*E_i=\sum_{j=1}^d \pi_j^*D_i$ which is ample on $X^d$.  Since $\phi$ is a finite surjective morphism, it follows that $E_i$ is ample.  By looking at the corresponding statement on $X^d$ we see that the intersection of any $dm+1$ distinct $E_i$'s is empty.  We also have $\dim \Sym X^d=dn$.  Since $r>2(dm)(dn)$, by Theorem \ref{cor4}(b) we have that all sets of $k$-rational $E$-integral points on $\Sym^d X$ are finite.  For $P\in R$ let $P^{(1)},\ldots,P^{(d)}$ denote the $d$ conjugates of $P$ over $k$.  Then $R'=\{(P^{(1)},\ldots,P^{(d)})\in X^d|P\in R\}$ is a set of $\sum_{i=1}^d\pi_i^*D$-integral points on $X^d$.  So $\phi(R')$ is a set of $E$-integral points on $\Sym^d X$.  Note that $\phi(R')$ is actually a set of $k$-rational points on $\Sym^d X$.  Therefore, from above, $\phi(R')$ must be finite, and so clearly $R$ must be finite.
\end{proof}
When $\#S=1$ we have the stronger theorem
\begin{theorem}
Let $X$ be a projective variety defined over $k=\mathbb{Q}$ or a complex quadratic field $k$.  Let $S=\{v_\infty\}$ consist of the unique archimedean place of $k$.  If $D_i$ is ample and defined over $k$ for all $i$ and $r>dm$ then all sets of $D$-integral points on $X$ of degree $d$ over $k$ are finite.
\end{theorem}
\begin{proof}
The proof is identical with the proof of the previous theorem, except that instead of using Theorem \ref{cor4}(b) we use Corollary \ref{qs1}.
\end{proof}
\section{A Result of Faltings}
\label{Faltings}
In \cite{Fa}, Faltings proves the finiteness of integral points on the complements of certain irreducible singular curves in $\mathbb{P}^2$.  Recently a similar result has also been obtained by Zannier in \cite{Co4}.  We show, as simple corollaries of our work on surfaces, how we may improve both results on integral points, and at the same time we will prove the analogous statement for holomorphic curves.  

Let $X$ be an irreducible nonsingular projective surface over an algebraically closed field $k$ of characteristic $0$.  Let $\ll=\co_X(L)$ be an ample line bundle on $X$ with $K_X+3L$ ample.  

Assume that the global sections $\Gamma(X,\ll)$ generate\\\\
(a).  $\ll_x/\mathfrak{m}_x^4\ll_x$ for all points $x\in X$\\
(b).  $\ll_x/\mathfrak{m}_x^3\ll_x \bigoplus \ll_y/\mathfrak{m}_y^3\ll_y$ for all pairs $\{x,y\}$ of distinct points\\
(c).  $\ll_x/\mathfrak{m}_x^2\ll_x \bigoplus \ll_y/\mathfrak{m}_y^2\ll_y \bigoplus \ll_z/\mathfrak{m}_z^2\ll_z$ for all triples $\{x,y,z\}$ of distinct points.\\

A three-dimensional subspace $E\subset \Gamma(X,\ll)$ that generates $\ll$ gives a morphism $f_E:X \to \mathbb{P}^2$.  Faltings studies this map when $E$ is suitably generic.
\begin{definition}
\label{Fgeneric}
Let $E\subset \Gamma(X,\ll)$ be a three-dimensional subspace.  We call $E$ generic if:\\\\
(a).  E generates $\ll$.\\
(b).  The discriminant locus $Z\subset X$ of $f_E$ is nonsingular.\\
(c).  The restriction of $f_E$ to $Z$ is birational onto its image $D\subset \mathbb{P}^2$.\\
(d).  $D$ has only cusps and nodes as singularities.
\end{definition}
Three-dimensional subspaces $E\subset \Gamma(X,\ll)$ are naturally parametrized by a Grassmannian $G$.  Let $n=L^2$.  It is then proven that
\begin{theorem}
\label{F2}
With notation as above\\
(a).  Generic $E$'s form a dense open subset $G'$ of $G$.\\
(b).  For generic $E$ let $\pi:Y\to X \to \mathbb{P}^2$ denote the associated normal Galois covering.  Then $Y$ is smooth, $Z$ is irreducible, and the covering group $Aut(Y/\mathbb{P}^2)$ is the full symmetric group $S_n$.
\end{theorem}

Faltings also proves
\begin{theorem}
\label{DP}
Let $\pi^*D$ be the pullback of $D$ to $Y$.  Then $\pi^*D=2\sum_{1\leq i<j \leq n}Z_{ij}=\sum_{i=1}^nA_i$ where $Z_{ij}$ is effective and nonsingular for every $i$ and $j$, and $A_i=\sum_{j\neq i}Z_{ij}$ is the pullback of $Z$ under the $i$th projection map $Y\to X$.  Furthermore, let $P\in \pi^*D$.  Then one of the following holds:\\\\
(a).  $\pi(P)$ is a smooth point of $D$ and $P\in Z_{ij}$ for exactly one $\{ij\}$.\\
(b).  $\pi(P)$ is a node of $D$ and exactly two components $Z_{ij}$ and $Z_{kl}$ of $\pi^*D$ for disjoint $\{ij\}$ and $\{kl\}$ intersect at $P$.\\
(c).  $\pi(P)$ is a cusp of $D$ and exactly three components $Z_{ij},Z_{ik},Z_{jk}$ intersect at $P$ for some $i,j,k$.
\end{theorem}
Let $d=\deg D$ and assume that everything above is defined over a number field.  The main result of \cite{Fa} is
\begin{theorem}[Faltings]
\label{Famain}
If $dL-\alpha Z$ is ample on $X$ for some $\alpha>12$ then $\mathbb{P}^2\backslash D$ is Mordellic.
\end{theorem}
Zannier proves this unconditionally if the Kodaira number of $X$ is nonnegative, and more generally he gives a numerical condition replacing the condition on $L$ and $Z$ above.  We will prove Theorem \ref{Famain} unconditionally, i.e. without the ampleness condition.  We also prove the analogue for holomorphic curves.  Under the assumptions discussed above, we prove
\renewcommand{\thetheorem}{\arabic{section}.\arabic{theorem}A}
\setcounter{theoremb}{\value{theorem}}
\begin{theorema}
$\mathbb{P}^2\backslash D$ is Mordellic.
\end{theorema}
\begin{theoremb}
$\mathbb{P}^2\backslash D$ is complete hyperbolic.  In particular, $\mathbb{P}^2\backslash D$ is Brody hyperbolic.
\end{theoremb}
\begin{proof}

Since $\pi:Y\backslash \pi^*D\to \mathbb{P}^2\backslash D$ is a finite \'etale covering, the problem is reduced to proving the theorems for $Y\backslash \pi^*D$.  The assumption (a) on $L$ given at the beginning of the section implies that $n=L^2\geq 9$.  We have $\pi^*D=\sum_{i=1}^n A_i$ and that $A_i$ is the pullback of $Z$ under the $i$th projection map $Y\to X$.  Therefore $A_i$ is ample as the projection is a finite map (recall that we assumed $Z\sim K_X+3L$ was ample).  It follows from Theorem \ref{DP} that at most four $A_i$'s meet at a point.  Therefore we're done by Theorems \ref{surf3a}(b) and \ref{surf3b}(b) with $r\geq 9$ and $m=4$.  That $\mathbb{P}^2\backslash D$ is complete hyperbolic follows from the fact that $Y\backslash \pi^*D$ is complete hyperbolic (see \cite{La2}).
\end{proof}

\section{Remarks on the Siegel and Picard-type Conjectures}
\renewcommand{\thetheorem}{\arabic{section}.\arabic{theorem}A}
\setcounter{theoremb}{\value{theorem}}
\label{Remarks}
In this section we will show the sharpness of the inequalities and the necessity of certain hypotheses in many of the conjectures, how our conjectures relate to other conjectures that have been made, and what special cases of the conjectures are known by previous work.
\subsection{Main Conjectures}
\subsubsection{Examples Limiting Improvements to the Conjectures}
Our main goal here is to show that the inequalities in all of the main conjectures cannot be improved.  We'll start with two fundamental examples on $\mathbb{P}^n$.
\begin{examplea}
\label{NHypera}
Let $X=\mathbb{P}^n$.  Let $D=\sum_{i=0}^nD_i$, where $D_i$ is the hyperplane defined by $x_i=0,i=0,\ldots,n$.  Let $k$ be a number field with an infinite number of units.  Let $S$ be the set of archimedean places.  Let $R$ be the set of points in $\mathbb{P}^n$ which have a representation where the coordinates are all units.  Then $R$ is a set of $D$-integral points on $X$.  It follows from the $S$-unit lemma that $R$ is Zariski-dense in $X$.
\end{examplea}
\begin{exampleb}
\label{NHyper}
Let $X$ and $D$ be as above.  Let $f_i,i=0,\ldots,n$ be linearly independent entire functions.  Let $f:\mathbb{C}\to X$ be defined by $f=(e^{f_0},\ldots,e^{f_n})$.  Clearly the image of $f$ does not intersect $D$.  It follows from Borel's lemma that the image of $f$ is Zariski-dense in $X$.
\end{exampleb}
We will give two variants of these examples which show that the inequalities in the Main Siegel and Picard-type Conjectures, Conjectures~\ref{conjmaina} and \ref{conjmainb}, are sharp for all values of $m$ and $\kappa_0$.
\begin{examplea}
\label{var1a}
Let $X$, $k$, $S$, $D$, $D_i$, and $R$ be as in Example \ref{NHypera}.  Let $Y=X^q$ and let $\pi_j$ be the $j$th projection map from $Y$ to $X$ for $j=1,\ldots,q$.  Let $R'=R^q\subset Y$.  Let $E_{i,j}=\pi_j^*D_i$ for $0\leq i \leq n,1\leq j \leq q$.  Let $1\leq m\leq nq$.  Let $r=\left[m+\frac{m}{n}\right]$ and $r'=\left[\frac{r}{n+1}\right]=\left[\frac{m}{n}\right]$.  Let 
\begin{equation*}
E=\sum_{j=1}^{r'}\sum_{i=0}^n E_{i,j}+\sum_{i=1}^{r-r'(n+1)}E_{i,r'+1}.
\end{equation*}
Then $R'$ is a set of $E$-integral points on $Y$ and it follows, again, from the $S$-unit lemma that $R'$ is Zariski-dense in $Y$.  Furthermore, there are at most $nr'+r-r'(n+1)=r-r'=m$ of the $E_{i,j}$'s in $E$ meeting at a given point, and $E$ is a sum of $r=\left[m+\frac{m}{n}\right]$ of the $E_{i,j}$'s with $\kappa(E_{i,j})=n$ for all $i$ and $j$.
\end{examplea}
\begin{exampleb}
\label{var1b}
Same as the above example, except that instead of $R'$, we use a holomorphic map $f:X\to Y\backslash E$ given by $f=(e^{f_{0,1}},\ldots,e^{f_{n,1}})\times\cdots \times(e^{f_{0,t}},\ldots,e^{f_{n,t}})$ where the $f_{i,j}$'s are linearly independent entire functions.  It follows from Borel's lemma that $f$ has Zariski-dense image in $Y$.
\end{exampleb}
The second variants of Examples \ref{NHypera} and \ref{NHyper} are
\begin{examplea}
\label{var2a}
Let $m$ and $n$ be positive integers.  Let $X$, $k$, $S$, $D_i$, $r$, and $r'$ be as in Example \ref{var1a}.  Let $D=\sum_{i=0}^{n}a_iD_i$ where $a_i=r'+1$ for $i=0,\ldots,r-(n+1)r'-1$ and $a_i=r'$ for $i=r-(n+1)r',\ldots,n$.  Then counting the $D_i$'s with their multiplicity in $D$, $D$ is a sum of $\sum_{i=0}^{n}a_i=r$ effective divisors such that the intersection of any $m+1$ of them is empty. We have $\kappa(D_i)=n$ for all $i$.  By Example \ref{NHypera} there exist dense sets of $D$-integral points on $X$.
\end{examplea}
\begin{exampleb}  
\label{var2b}
The same example as above, except we use the holomorphic map from Example \ref{NHyper}.
\end{exampleb}
The above four examples also show that one cannot improve the inequalities in Conjectures \ref{conj1a},B and \ref{conj1ab},B.

We have not yet discussed the $\kappa_0=0$ case.  If $D$ is a divisor on a projective variety $X$, then by blowing up subvarieties of $D$ on $X$ we may get a divisor $D'$ on $X'$ with arbitrarily many components and $X\backslash D \cong X'\backslash D'$.  In this case, the new components $C$ have $\kappa(C)=0$.  So, as is suggested by the $\kappa_0$ in the denominators of the inequalities, there is no possible result of the type in the Main Siegel and Picard-type Conjectures if one allows divisors $D_i$ with $\kappa(D_i)=0$.  However, all is not lost in this case.  If we are willing to include in the inequalities numerical invariants of the variety such as the Picard number, then it is possible to give theorems for arbitrary effective divisors.  We will discuss this in Section \ref{mainknown}.

There are also examples showing that the exceptional sets may be dense, even if the hypotheses of the Main Siegel and Picard-type Conjectures are satisfied.  For example, let $X=\mathbb{P}^1\times \mathbb{P}^1$ and let $D=\sum_{i\in I} P_i\times \mathbb{P}^1$ be a finite sum with $P_i\in \mathbb{P}^1(k), i \in I$, for some number field $k$.  Then it is easy to show that $\Excd(X\backslash D)=\Exch(X\backslash D)=X\backslash D$.

For the Main Conjectures for Ample Divisors we have
\begin{examplea}
Let $D$ be the sum of any $r$ hyperplanes in general position (i.e. the intersection of any $n+1$ of them is empty) in $\mathbb{P}^n$ with $n<r\leq 2n$.  Assume also that $D$ is defined over a number field.  Then one may show that there exists a linear subspace $L\subset \mathbb{P}^n$ with $\dim L=\left[\frac{n}{r-n}\right]$ such that $L$ contains a dense set of $D|_L$-integral points (for some $k$ and $S$) (see \cite{Fu2}, \cite{Gr}, and \cite{No} for the constructions).
\end{examplea}
\begin{exampleb}
In the same situation as above, one may also show that there exists a holomorphic map $f:\mathbb{C}\to L\backslash D$ with Zariski-dense image.
\end{exampleb}
In the simplest case, where $r=2m=2n$, we may simply take $L$ to be a line that passes through points $P$ and $Q$ where $P$ is the intersection of, say, the first $n$ hyperplanes and $Q$ is the intersection of the last $n$ hyperplanes.  Then $L\cap D$ is a $\mathbb{P}^1$ minus two points, and so we see that we cannot have finiteness or constancy for the objects in question.
\renewcommand{\thetheorem}{\arabic{section}.\arabic{theorem}}
\begin{remark}
\label{rbig}
It is quite possible that our Main Conjectures for Ample Divisors may be extended to quasi-ample divisors.  Let $D$ be a quasi-ample divisor on a projective variety $X$.  Let $n>0$ be large enough such that the map $\Phi=\Phi_{nD}$, corresponding to $nD$, is birational.  It is then quite plausible that all of our conclusions that held for ample divisors generalize to quasi-ample divisors if we state things in terms of $\Phi$, that is, replace $\dim \Excd(X\backslash D)$ and $\dim \Exch(X\backslash D)$ by $\dim \Phi(\Excd(X\backslash D))$ and $\dim \Phi(\Excd(X\backslash D))$ in the conjectures.
\end{remark}

\subsubsection{Relation to Vojta's Main Conjecture}
\label{mainrelation}
We now show how some special cases of the Main Conjectures are related to Vojta's Main Conjecture.  If $D$ is a divisor on a nonsingular complex variety $X$, we say that $D$ has normal crossings if every point $P\in D$ has an analytic open neighborhood in $X$ with analytic local coordinates $z_1,\ldots,z_n$ such that $D$ is locally defined by $z_1\cdot z_2\cdots z_i=0$ for some $i$.  Inspired by results in equi-dimensional Nevanlinna theory, Vojta made the following conjecture in \cite{Vo2}.
\renewcommand{\thetheorem}{\arabic{section}.\arabic{theorem}A}
\setcounter{theoremb}{\value{theorem}}
\begin{conjecturea}[Vojta's Main Conjecture]
\label{Vmain}
Let $X$ be a nonsingular projective variety with canonical divisor $K$.  Let $D$ be a normal crossings divisor on $X$, and let $k$ be a number field over which $X$ and $D$ are defined.  Let $A$ be a quasi-ample divisor on $X$.  Let $\epsilon>0$.  Then there exists a proper Zariski-closed subset $Z=Z(X,D,\epsilon,A)$ such that
\begin{equation*}
m(D,P)+h_K(P)\leq \epsilon h_A(P)+O(1)
\end{equation*}
for all points $P\in X\backslash Z$.
\end{conjecturea}
Similarly, the analogue is conjectured for holomorphic curves
\begin{conjectureb}
\label{Vmainb}
Let $X$ be a nonsingular complex projective variety with canonical divisor $K$.  Let $D$ be a normal crossings divisor on $X$.  Let $A$ be a quasi-ample divisor on $X$.  Let $\epsilon>0$.  Then there exists a proper Zariski-closed subset $Z=Z(X,D,\epsilon,A)$ such that for all holomorphic maps $f:\mathbb{C}\to X$ whose image is not contained in $Z$,
\begin{equation*}
m(D,r)+T_K(r)\leq \epsilon T_A(r)+O(1)
\end{equation*}
holds for all $r$ outside a set of finite Lebesgue measure.
\end{conjectureb}
Qualitatively, these conjectures have the following simple consequences.
\begin{conjecturea}
\label{conj3}
Let $X$ be a nonsingular projective variety, defined over a number field $k$.  Let $K$ be the canonical divisor of $X$, and $D$ a normal crossings divisor on $X$, defined over $k$.  Suppose that $K+D$ is quasi-ample. Then $X\backslash D$ is quasi-Mordellic.
\end{conjecturea}
\begin{conjectureb}
\label{conj3b}
Let $X$ be a nonsingular complex projective variety.  Let $K$ be the canonical divisor of $X$, and $D$ a normal crossings divisor on $X$.  Suppose that $K+D$ is quasi-ample.  Then $X\backslash D$ is quasi-Brody hyperbolic.
\end{conjectureb}
\renewcommand{\thetheorem}{\arabic{section}.\arabic{theorem}}
To relate these conjectures to our conjectures we recall the following theorem, which is a consequence of Mori theory \cite[Lemma 1.7]{Mo}.
\begin{theorem}
Let $X$ be a nonsingular complex projective variety of dimension $n$.  If $D_1,\ldots,D_{n+2}$ are ample divisors on $X$  then $K+\sum_{i=1}^{n+2}D_i$ is ample.
\end{theorem}
So when $X$ is nonsingular, the $D_i$'s are ample, and $D$ has normal crossings, we see that Conjectures \ref{conj1ab} and \ref{conj1bb} are consequences of Conjectures \ref{conj3} and \ref{conj3b}.

\subsubsection{Previously Known Results Related to the Conjectures}
\label{mainknown}
As was discussed earlier, our work builds on previous work of Corvaja and Zannier, who obtained results on surfaces in \cite{Co2}, and initiated the general method we have used in \cite{Co}.  The Nevanlinna theoretic analogues of \cite{Co2} were proved by Liu and Ru in \cite{Ru4}.  We briefly discussed these previous results in Section \ref{ssurf}.

We now discuss what is known for arbitrary divisors.  As a consequence of his work on integral points on subvarieties of semi-abelian varieties, Vojta \cite{Vo1} proved
\renewcommand{\thetheorem}{\arabic{section}.\arabic{theorem}A}
\setcounter{theoremb}{\value{theorem}}
\begin{theorema}
\label{Vojtaa}
Let $X$ be a projective variety defined over a number field $k$.  Let $\rho$ denote the Picard number of $X$.  Let $D$ be an effective divisor on $X$ defined over $k$ which has more than $\dim X - h^1(X,\mathcal{O}_X)+\rho$ (geometrically) irreducible components.  Then $X\backslash D$ is quasi-Mordellic.
\end{theorema}
Similarly, a special case of work of Noguchi \cite{No2} gives
\begin{theoremb}
\label{Vojtab}
Let $X$ be a complex projective variety.  Let $\rho$ denote the Picard number of $X$.  Let $D$ be an effective divisor on $X$ which has more than $\dim X - h^1(X,\mathcal{O}_X)+\rho$ irreducible components.  Then $X\backslash D$ is quasi-Brody hyperbolic.
\end{theoremb}
We note that it is easily shown that both theorems are sharp in that there are divisors with $\dim X - h^1(X,\mathcal{O}_X)+\rho$ irreducible components for which the conclusions of the theorems are false.  For a weaker, but more elementary theorem along these lines, see also Th. 2.4.1 in~\cite{Vo2}.  As consequences of Theorems \ref{Vojtaa},B we see that Conjectures \ref{conj1ab},B are true for $X=\mathbb{P}^n$, and more generally for any projective variety $X$ with Picard number one.

From the work of Noguchi and Winkelmann \cite{No} we have the following theorems related to our Main Conjectures for Ample Divisors (some special cases of these results had been obtained previously by various people; see \cite{No} for the history).
\begin{theorema}
Let $X$ be a projective variety of dimension $n$ defined over a number field $k$.  Let $S$ be a finite set of places of $k$ containing the archimedean places.  Let $\rho$ be the Picard number of $X$.  Let $D=\sum_{i=1}^rD_i$ be a divisor on $X$ defined over $k$ with the $D_i$'s effective reduced ample Cartier divisors such that the intersection of any $n+1$ of them is empty.\\\\
(a).  If $r>n+1$ then all sets of $D$-integral points $R$ have $\dim R\leq \frac{n}{r-n}\rho$.\\
(b).  If $r>n(\rho+1)$ then $X\backslash D$ is Mordellic.\\
(c).  If $X\subset \mathbb{P}^N$, all $D_i$'s are hypersurface cuts of $X$, and $r>2n$ then $X\backslash D$ is Mordellic.
\end{theorema}
\begin{theoremb}
Let $X$ be a complex projective variety of dimension $n$.  Let $\rho$ be the Picard number of $X$.  Let $D=\sum_{i=1}^rD_i$ be a divisor on $X$ with the $D_i$'s effective reduced ample Cartier divisors such that the intersection of any $n+1$ of them is empty.\\\\
(a).  If $r>n+1$ then all holomorphic maps $f:\mathbb{C}\to X\backslash D$ have $\dim f(\mathbb{C})\leq \frac{n}{r-n}\rho$.\\
(b).  If $r>n(\rho+1)$ then $X\backslash D$ is complete hyperbolic and hyperbolically imbedded in $X$.  In particular, $X\backslash D$ is Brody hyperbolic.\\
(c).  If $X\subset \mathbb{P}^N$, all $D_i$'s are hypersurface cuts of $X$, and $r>2n$ then $X\backslash D$ is complete hyperbolic and hyperbolically imbedded in $X$.  In particular, $X\backslash D$ is Brody hyperbolic.
\end{theoremb}
Consequently, when $m=\dim X$, the $D_i$'s are reduced divisors, and $\rho(X)=1$, we have that the Main Conjectures for Ample Divisors, Conjectures \ref{conj2a},B, are true modulo the statements on the exceptional sets (i.e. replace $\Excd(X\backslash D)$ by any particular set of integral points $R$ in Conjecture \ref{conj2a}, etc.)  Similarly, the part (c)'s of the above theorems give special cases of the part (b)'s of Conjectures \ref{conj2a},B.
\subsection{General Conjectures}
\renewcommand{\thetheorem}{\arabic{section}.\arabic{theorem}}
\subsubsection{Examples Limiting Improvements to the Conjectures}
We start off with an example showing that the inequalities in the General Conjectures are best possible when $X$ is a curve.
\begin{example}
\label{genex}
Let $X$ be a projective curve defined over a number field $k$ with $\co_k^*$ infinite.  Let $f:X\to \mathbb{P}^1$ be a morphism of degree $d$ defined over $k$.  Let $P,Q\in \mathbb{P}^1(k)$ be two distinct points over which $f$ is unramified, and let $D=P+Q$.  Then there exists an infinite set $R$ of $k$-rational $D$-integral points on $\mathbb{P}^1\backslash D$.  Since $f$ has degree $d$, $f^{-1}(R)$ is a set of $f^*D$-integral points on $X\backslash f^*D$ of degree $d$ over $k$ and $f^*D$ is a sum of $2d$ distinct points on $X$.
\end{example}
Taking products of curves, we then get examples in all dimensions showing that the inequality in the General Siegel-type Conjecture cannot be improved in the case $\kappa_0=1$.
\begin{example}
\label{genex2}
Let $D=\sum_{i=1}^{2md}H_i$ be a sum of hyperplanes on $\mathbb{P}^n$ defined over a number field $k$ such that the intersection of any $m+1$ of the $H_i$'s is empty.  Suppose also that $\bigcap_{i=(j-1)m+1}^{jm}H_i=\{P_j\}$ consists of a single point for $j=1,\ldots,2d$ and the $P_j$'s are collinear.  Then there exist infinite sets of $D$-integral points of degree $d$ on $\mathbb{P}^n\backslash D$ over large enough number fields.  Indeed, the line $L$ through the $P_j$'s intersects $D$ in $2d$ points, and we see from Example \ref{genex} that $L\backslash L\cap D$ contains infinite sets of integral points over large enough number fields.
\end{example}
This shows that the inequality in the finiteness part of the General Siegel-type Conjecture for Ample Divisors cannot be improved.  We expect that using only divisors that are sums of hyperplanes on projective space, one may show that the other inequalities in the General Conjectures may not be improved for any set of parameters.  For example, it should be true that if $D$ is a sum of $2d+n-1$ hyperplanes in general position on $\mathbb{P}^n$, then for some number field $k$ there exist dense sets of $D$-integral points on $\mathbb{P}^n$ of degree $d$ over $k$.  In any case, it is easy to show that $\Excd(\mathbb{P}^n\backslash D)=\mathbb{P}^n\backslash D$.  If $P$ is a point where $n$ of the hyperplanes intersect, then any line through $P$ will intersect $D$ in $2d$ points.  But as we have seen, over some number field $k$, such lines will contain infinitely many integral points of degree $d$ over $k$.  To show the existence of a Zarisk-dense set of $D$-integral points, one needs to show that if the lines and their sets of integral points are chosen correctly, then the infinite union of the sets of integral points will still be a set of $D$-integral points (there is no problem for finite unions).

\subsubsection{Vojta's General Conjecture and a Conjectural Discriminant-Height Inequality}
We will now investigate how the General Siegel-type Conjecture, Conjecture \ref{congen}, is related to Vojta's General Conjecture.  In order to make a connection between the two conjectures, we will need to formulate a new conjecture bounding the absolute logarithmic discriminant in terms of heights.  We will digress briefly to discuss this new conjecture.  Let $X$ be a variety defined over a number field $k$ and let $P\in X(\kk)$.  Let $d(P)=\frac{1}{[k(P):\mathbb{Q}]}\log |D_{k(P)/\mathbb{Q}}|$ where $D_{k(P)/\mathbb{Q}}$ is the discriminant of $k(P)$ over $\mathbb{Q}$.  We call $d(P)$ the absolute logarithmic discriminant of $P$.  Let $m(D,P)=\sum_{v\in S}\lambda_{D,v}(P)$.  Then Vojta's General Conjecture states
\begin{conjecture}[Vojta's General Conjecture]
\label{Vgeneral}
Let $X$ be a complete nonsingular variety with canonical divisor $K$.  Let $D$ be a normal crossings divisor on $X$, and let $k$ be a number field over which $X$ and $D$ are defined.  Let $A$ be a quasi-ample divisor on $X$.  Let $\epsilon>0$.  If $\nu$ is a positive integer then there exists a Zariski-closed subset $Z=Z(\nu,X,D,\epsilon,A)$ such that
\begin{equation*}
m(D,P)+h_K(P)\leq d(P)+\epsilon h_A(P)+O(1)
\end{equation*}
for all points $P\in X(\kk)\backslash Z$ such that $[k(P):k]\leq \nu$.
\end{conjecture}
Actually, Vojta's General Conjecture as it appears in \cite{Vo2} has the discriminant term as $(\dim X)d(P)$, but it is now believed that the $\dim X$ term is unecessary (see \cite[Conjecture 8.7]{Vo5} or the discussion at the end of \cite{Vo6}).

Vojta's General Conjecture, with $D=0$, can be seen as giving a lower bound on the absolute logarithmic discriminant in terms of heights (outside some Zariski-closed subset).  As a companion to this, we give the following conjectural upper bound on the logarithmic discriminant in terms of heights.
\begin{conjecture}
\label{conj4}
Let $X$ be a nonsingular projective variety of dimension $n$ defined over a number field $k$ with canonical divisor $K$.  Let $A$ be an ample divisor on $X$.  Let $\nu$ be a positive integer.  Let $\epsilon>0$.  Then
\begin{equation*}
d(P)\leq h_K(P)+(2[k(P):k]+n-1+\epsilon)h_A(P)+O(1)
\end{equation*}
for all $P\in X(\kk)$ with $[k(P):k]\leq \nu$.
\end{conjecture}
\begin{remark}
It is possible that with the hypothesis $A$ ample weakened to $A$ quasi-ample that the inequality holds outside of some Zariski-closed subset of $X$ (it is not hard to see the necessity of the Zariski-closed subset in this case).  It is also possible that the conjecture is true with $\epsilon=0$.  As with Vojta's General Conjecture, it is quite plausible that one may take $\nu=\infty$, i.e. the inequality holds for all $P\in X(\kk)$.
\end{remark}
It is a result of Silverman \cite{Si2} that Conjecture \ref{conj4} is true for $X=\mathbb{P}^n$ with $\epsilon=0$ and $\nu=\infty$.  For curves, Conjecture \ref{conj4} is true by a result of Song and Tucker \cite[Eq. 2.0.3]{Tu}.  They proved the stronger statement
\begin{theorem}
\label{thTu}
Let $X$ be a nonsingular projective curve defined over a number field $k$ with canonical divisor $K$.  Let $A$ be an ample divisor on $X$.  Let $\nu$ be a positive integer.  Let $\epsilon>0$.  Then
\begin{equation*}
d(P)\leq d_a(P)\leq h_K(P)+(2[k(P):k]+\epsilon)h_A(P)+O(1)
\end{equation*}
for all $P\in X(\kk)$ with $[k(P):k]\leq \nu$, where $d_a(P)$ is the arithmetic discriminant of $P$ (see \cite{Vo8} for the definition and properties).
\end{theorem}
We now show how Vojta's General Conjecture, combined with our conjectural upper bound on the discrimant, imply a special case of the General Siegel-type Conjecture.
\begin{theorem}
Assume Vojta's General Conjecture, Conjecture \ref{Vgeneral}, and the conjectural upper bound on the absolute logarithmic discriminant, Conjecture \ref{conj4}.  Let $X$ be a nonsingular projective variety defined over a number field $k$.  Let $n=\dim X$.  Let $D=\sum_{i=1}^r D_i$ be a normal crossings divisor defined over $k$ with $D_i$ ample and effective for all $i$.  If $r>2\nu+n-1$ then $X\backslash D$ is degree $\nu$ quasi-Mordellic.  In particular, there do not exist Zariski-dense sets of $D$-integral points on $X$ of degree $\nu$ over $k$.
\end{theorem}
\begin{proof}
Let $R$ be a set of $D$-integral points on $X$ of degree $\nu$ over $k$.  Then $m(D,P)+h_K(P)=h_D(P)+h_K(P)+O(1)$ for $P\in R$.  By Conjecture \ref{conj4}, for any $\epsilon>0$, $h_{D_i}(P)\geq \frac{d(P)-h_K(P)}{2\nu+n-1+\epsilon}+O(1)$.  So since $r>2\nu+n-1$, we have $h_D(P)\geq d(P)-h_K(P)+(1-\epsilon)h_{D_1}(P)+O(1)$.  Therefore 
\begin{equation*}
m(D,P)+h_K(P)>d(P)+\epsilon h_A(P)+O(1)
\end{equation*}
for all $P\in R$, for any ample divisor $A$ on $X$ and small enough $\epsilon$.  So we're done by Vojta's General Conjecture.
\end{proof}
So assuming Vojta's General Conjecture and Conjecture \ref{conj4}, we see that the General Siegel-type Conjecture is true if $D_i$ is ample for all $i$ and $D$ has normal crossings.

\subsubsection{Previously Known Results Related to the Conjectures}
\label{sVo}
In \cite{Vo7}, Vojta proved the following generalization of Falting's theorem on rational points on curves and the Thue-Siegel-Roth-Wirsing theorem.
\begin{theorem}
Let $X$ be a nonsingular projective curve defined over a number field $k$ with canonical divisor $K$.  Let $D$ be an effective divisor on $X$ defined over $k$ with no multiple components and $A$ an ample divisor on $X$.  Let $\nu$ be a positive integer and let $\epsilon>0$.  Then
\begin{equation*}
m(D,P)+h_K(P)\leq d_a(P)+\epsilon h_A(P)+O(1)
\end{equation*}
for all $P\in X(\kk)\backslash D$ with $[k(P):k]\leq \nu$, where the constant in $O(1)$ depends on $X,D,\nu,A$, and $\epsilon$.
\end{theorem}
Using Theorem \ref{thTu} we then easily obtain the following theorem.
\begin{corollary}
Let $X$ be a nonsingular projective curve defined over a number field $k$.  Let $D$ be an effective divisor on $X$ that is a  sum of more than $2\nu$ distinct points.  Then $X\backslash D$ is degree $\nu$ Mordellic.
\end{corollary}
Therefore our General Siegel-type Conjectures are true for curves.  Of course for $\mathbb{P}^1$ this was already known from the Thue-Siegel-Roth-Wirsing theorem.  As mentioned earlier, the special case $\nu=2$ was also proven by Corvaja and Zannier using the Schmidt Subspace Theorem technique \cite{Co2}.

\subsection{Conjectures over $\mathbb{Z}$ and Complex Quadratic Rings of Integers}
I am not aware of any previous results that pertain to these conjectures, or any way to relate them to other known conjectures.  An open problem then is to formulate quantitative conjectures explaining the qualitative conjectures I have made over $\mathbb{Z}$ and complex quadratic rings of integers.  We now briefly discuss some examples showing that in many cases the inequalities in these conjectures may not be improved.

For the Main Siegel-type Conjecture over $\mathbb{Z}$, to show that the inequality in the conjecture may not be improved we may simple take $D=mH$ where $H$ is a hyperplane on $\mathbb{P}^n$ defined over $\mathbb{Q}$.  Examples where the $m$ divisors have no components in common are easily obtained from products of projective spaces.  

For the Main Conjecture on Ample Divisors over $\mathbb{Z}$, if $D=\sum_{i=1}^mH_i$ is a sum of $m<n$ distinct hyperplanes on $\mathbb{P}^n$ defined over $\mathbb{Q}$ then $\dim \cap_{i=1}^mH_i=n-m$ and there is a $Y=\mathbb{P}^{n-m+1}\subset \mathbb{P}^n$ with $D|_Y$ a hyperplane on $Y$ defined over $\mathbb{Q}$.  So there are sets of $D$-integral points on $\mathbb{P}^n$ with dimension $n-m+1$.

Examples for the General Conjectures over $\mathbb{Z}$ are nearly identical to Examples \ref{genex} and \ref{genex2}, except that we must replace $2d$ by $d$ everywhere, since we are using $\mathbb{A}^1$ as our starting point.  Again, we expect that using only divisors that are sums of hyperplanes on projective space, one may show that the inequalities in the General Conjectures over $\mathbb{Z}$ may not be improved for any set of parameters.
\bibliography{integral}

\begin{thebibliography}{10}

\bibitem{Co6}
Pietro Corvaja and Umberto Zannier.
\newblock On the integral points on certain surfaces.
\newblock {\em To appear}.

\bibitem{Co}
Pietro Corvaja and Umberto Zannier.
\newblock A subspace theorem approach to integral points on curves.
\newblock {\em C. R. Math. Acad. Sci. Paris}, 334(4):267--271, 2002.

\bibitem{Co5}
Pietro Corvaja and Umberto Zannier.
\newblock On the number of integral points on algebraic curves.
\newblock {\em J. Reine Angew. Math.}, 565:27--42, 2003.

\bibitem{Co3}
Pietro Corvaja and Umberto Zannier.
\newblock On a general {T}hue's equation.
\newblock {\em Amer. J. Math.}, 126(5):1033--1055, 2004.

\bibitem{Co2}
Pietro Corvaja and Umberto Zannier.
\newblock On integral points on surfaces.
\newblock {\em Ann. of Math.}, 160(2):705--726, 2004.

\bibitem{Fa}
Gerd Faltings.
\newblock A new application of {D}iophantine approximations.
\newblock In {\em A panorama of number theory or the view from Baker's garden
  (Z\"urich, 1999)}, pages 231--246. Cambridge Univ. Press, Cambridge, 2002.

\bibitem{Fa2}
Gerd Faltings and Gisbert W{\"u}stholz.
\newblock Diophantine approximations on projective spaces.
\newblock {\em Invent. Math.}, 116(1-3):109--138, 1994.

\bibitem{Fu2}
Hirotaka Fujimoto.
\newblock Extensions of the big {P}icard's theorem.
\newblock {\em T\^ohoku Math. J. (2)}, 24:415--422, 1972.

\bibitem{Gr}
Mark~L. Green.
\newblock Holomorphic maps into complex projective space omitting hyperplanes.
\newblock {\em Trans. Amer. Math. Soc.}, 169:89--103, 1972.

\bibitem{Gr2}
Mark~L. Green.
\newblock The hyperbolicity of the complement of {$2n+1$} hyperplanes in
  general position in {$P\sb{n}$} and related results.
\newblock {\em Proc. Amer. Math. Soc.}, 66(1):109--113, 1977.

\bibitem{Ii}
Shigeru Iitaka.
\newblock {\em Algebraic geometry}, volume~76 of {\em Graduate Texts in
  Mathematics}.
\newblock Springer-Verlag, New York, 1982.
\newblock An introduction to birational geometry of algebraic varieties,
  North-Holland Mathematical Library, 24.

\bibitem{Ka}
Yujiro Kawamata.
\newblock A generalization of {K}odaira-{R}amanujam's vanishing theorem.
\newblock {\em Math. Ann.}, 261(1):43--46, 1982.

\bibitem{Kl}
Steven~L. Kleiman.
\newblock Toward a numerical theory of ampleness.
\newblock {\em Ann. of Math. (2)}, 84:293--344, 1966.

\bibitem{La}
Serge Lang.
\newblock {\em Fundamentals of {D}iophantine geometry}.
\newblock Springer-Verlag, New York, 1983.

\bibitem{La2}
Serge Lang.
\newblock {\em Introduction to complex hyperbolic spaces}.
\newblock Springer-Verlag, New York, 1987.

\bibitem{Ru4}
Yuancheng Liu and Min Ru.
\newblock Degeneracy of holomorphic curves in surfaces.
\newblock {\em Preprint}.

\bibitem{Mo}
Shigefumi Mori.
\newblock Threefolds whose canonical bundles are not numerically effective.
\newblock {\em Ann. of Math. (2)}, 116(1):133--176, 1982.

\bibitem{No2}
Junjiro Noguchi.
\newblock Lemma on logarithmic derivatives and holomorphic curves in algebraic
  varieties.
\newblock {\em Nagoya Math. J.}, 83:213--233, 1981.

\bibitem{No}
Junjiro Noguchi and J{\"o}rg Winkelmann.
\newblock Holomorphic curves and integral points off divisors.
\newblock {\em Math. Z.}, 239(3):593--610, 2002.

\bibitem{Ru}
Min Ru.
\newblock {\em Nevanlinna theory and its relation to {D}iophantine
  approximation}.
\newblock World Scientific Publishing Co. Inc., River Edge, NJ, 2001.

\bibitem{Ru3}
Min Ru.
\newblock A defect relation for holomorphic curves intersecting hypersurfaces.
\newblock {\em Amer. J. Math.}, 126(1):215--226, 2004.

\bibitem{Si2}
Joseph~H. Silverman.
\newblock Lower bounds for height functions.
\newblock {\em Duke Math. J.}, 51(2):395--403, 1984.

\bibitem{Tu}
Xiangjun Song and Thomas~J. Tucker.
\newblock Dirichlet's theorem, {V}ojta's inequality, and {V}ojta's conjecture.
\newblock {\em Compositio Math.}, 116(2):219--238, 1999.

\bibitem{Vo2}
Paul Vojta.
\newblock {\em Diophantine approximations and value distribution theory},
  volume 1239 of {\em Lecture Notes in Mathematics}.
\newblock Springer-Verlag, Berlin, 1987.

\bibitem{Vo6}
Paul Vojta.
\newblock A refinement of {S}chmidt's subspace theorem.
\newblock {\em Amer. J. Math.}, 111(3):489--518, 1989.

\bibitem{Vo8}
Paul Vojta.
\newblock Arithmetic discriminants and quadratic points on curves.
\newblock In {\em Arithmetic algebraic geometry (Texel, 1989)}, volume~89 of
  {\em Progr. Math.}, pages 359--376. Birkh\"auser Boston, Boston, MA, 1991.

\bibitem{Vo7}
Paul Vojta.
\newblock A generalization of theorems of {F}altings and
  {T}hue-{S}iegel-{R}oth-{W}irsing.
\newblock {\em J. Amer. Math. Soc.}, 5(4):763--804, 1992.

\bibitem{Vo1}
Paul Vojta.
\newblock Integral points on subvarieties of semiabelian varieties. {I}.
\newblock {\em Invent. Math.}, 126(1):133--181, 1996.

\bibitem{Vo3}
Paul Vojta.
\newblock On {C}artan's theorem and {C}artan's conjecture.
\newblock {\em Amer. J. Math.}, 119(1):1--17, 1997.

\bibitem{Vo5}
Paul Vojta.
\newblock Nevanlinna theory and {D}iophantine approximation.
\newblock In {\em Several complex variables (Berkeley, CA, 1995--1996)},
  volume~37 of {\em Math. Sci. Res. Inst. Publ.}, pages 535--564. Cambridge
  Univ. Press, Cambridge, 1999.

\bibitem{Co4}
Umberto Zannier.
\newblock On the integral points on the complement of ramification-divisors (to
  appear).
\newblock {\em Jussieu Math. J.}

\end{thebibliography}
\end{document}